\documentclass{article}
\usepackage{appendix}
\usepackage{mathrsfs}
\usepackage{amsmath}
\usepackage{amsthm}
\usepackage{amssymb}
\usepackage{graphicx}
\usepackage{fancyhdr}
\usepackage [latin1]{inputenc}
\usepackage{enumerate}
\usepackage{geometry}
\catcode`\@=11 \@addtoreset{equation}{section}

\catcode`\@=12
\usepackage{color}

\usepackage[colorlinks,
            linkcolor=blue,
            anchorcolor=blue,
            citecolor=red,
            ]{hyperref}

\DeclareFontFamily{U}{mathx}{}
\DeclareFontShape{U}{mathx}{m}{n}{ <-> mathx10 }{}
\DeclareSymbolFont{mathx}{U}{mathx}{m}{n}
\DeclareFontSubstitution{U}{mathx}{m}{n}

\DeclareMathAccent{\widecheck}{0}{mathx}{"71}

\makeatletter
\newcommand{\wcheck}[1]{\mathpalette\wcheck@{#1}}
\newcommand{\wcheck@}[2]{%
  \begingroup
  \edef\wcheck@font{\the
    \ifx#1\displaystyle\textfont\else\ifx#1\textstyle\textfont
    \else\ifx#1\scriptstyle\scriptfont\else\scriptscriptfont\fi\fi\fi\@ne
  }%
  \sbox\z@{\wcheck@font\mbox{#2}\mbox{\char\the\skewchar\font}}%
  \sbox\tw@{\wcheck@font#2\char\the\skewchar\font}%
  \dimen@=\dimexpr\wd\tw@-\wd\z@\relax
  {\,\kern2\dimen@\widecheck{\!\kern-2\dimen@#2\!}\,}%
  \endgroup
}
\makeatother

\newtheorem{Theorem}{Theorem}[section]
\newtheorem{Proposition}{Proposition}[section]
\newtheorem{Lemma}{Lemma}[section]
\newtheorem{Corollary}{Corollary}[section]

\newcommand{\newcom}{\newcommand}
\newcommand{\bTheorem}[1]{
\begin{Theorem} \label{T#1} }
\newcommand{\eT}{\end{Theorem}}

\newcommand{\bProposition}[1]{
\begin{Proposition} \label{P#1}}
\newcommand{\eP}{\end{Proposition}}

\newcommand{\bLemma}[1]{
\begin{Lemma} \label{L#1} }
\newcommand{\eL}{\end{Lemma}}

\newcommand{\bCorollary}[1]{
\begin{Corollary} \label{C#1} }
\newcommand{\eC}{\end{Corollary}}

\newcommand{\beq}{\begin{equation}}
\newcommand{\eeq}{\end{equation}}
\newcom{\ben}{\begin{eqnarray}}
\newcom{\een}{\end{eqnarray}}
\newcom{\beno}{\begin{eqnarray*}}
\newcom{\eeno}{\end{eqnarray*}}
\newcom{\bali}{\begin{aligned}}
\newcom{\eali}{\end{aligned}}

\newcommand{\bFormula}[1]{
\begin{equation} \label{#1}}
\newcommand{\eF}{\end{equation}}

\newcommand{\Grad}{\nabla}

\newcommand{\ep}{\varepsilon}

\newcommand{\norm}[1]{{\left\|#1\right\|}}

\newcommand{\dd}{\,{\rm d}}

\newcommand\Cbox[2]{%
    \newbox\contentbox%
    \newbox\bkgdbox%
    \setbox\contentbox\hbox to \hsize{%
        \vtop{
            \kern\columnsep
            \hbox to \hsize{%
                \kern\columnsep%
                \advance\hsize by -2\columnsep%
                \setlength{\textwidth}{\hsize}%
                \vbox{
                    \parskip=\baselineskip
                    \parindent=0bp
                    #2
                }%
                \kern\columnsep%
            }%
            \kern\columnsep%
        }%
    }%
    \setbox\bkgdbox\vbox{
        \color{#1}
        \hrule width  \wd\contentbox %
               height \ht\contentbox %
               depth  \dp\contentbox
        \color{black}
    }%
    \wd\bkgdbox=0bp%
    \vbox{\hbox to \hsize{\box\bkgdbox\box\contentbox}}%
    \vskip\baselineskip%
}

\usepackage{graphicx} % Required for inserting images

\title{\bf{Sharp asymptotic stability of the incompressible porous medium equation on the whole plane }}

\author{
Sai Li\thanks{Department of Applied Mathematics, Nanjing Forestry University, Nanjing--210037, People's Republic of China, E-mail: \texttt{lsmath@njfu.edu.cn}}
\qquad  \qquad
Yang Li\thanks{School of Mathematical Sciences, Anhui University, Hefei--230601, People's Republic of China, E-mail: \texttt{lynjum@163.com}} 
}

\begin{document}
\maketitle 
\abstract{
We study the asymptotic stability of the steady state
$e_s(x)=-x_2$ of the incompressible porous media equation
in $\mathbb R^2$.
For every $s>2$, we prove that sufficiently small initial
perturbations $\theta_0\in H^s(\mathbb R^2)$ of $e_s$
generate unique global solutions $\rho=e_s+\theta$,
where $\theta$ remains uniformly bounded in $H^s$
and converges strongly to zero in $H^s$ as $t\to\infty$.
The main innovations of this work are threefold.
First, exploiting the structure of the equation,
we construct a weighted energy that eliminates the need
to control $\|\Grad u^2\|_{L^1_tL^\infty_x}$.
Second, rather than estimating the anisotropic commutators
directly, we study the integral terms containing these
commutators and establish suitable bounds through
refined computations in frequency space.
Third, we establish convergence in the highest-order
$H^s$ norm under uniform energy bounds involving only
an $L^\infty_tH^s_x$ bound on $\theta$
and an $L^2_tH^s_x$ bound on $\mathcal R_1\theta$.
%To the best of our knowledge, this is the first global
%stability result in $\mathbb R^2$ covering small
%$H^s$ perturbations for every $s>2$.
}

\medskip

\vspace{2mm}


\tableofcontents

\section{Introduction}

The two-dimensional incompressible porous media (IPM) system
describes the transport of fluid density through a porous medium
under gravity and takes the form
\begin{equation}
\label{equ1.1}
\begin{cases}
\partial_t \rho
+
u
\cdot
\nabla \rho
=
0,
\\[2mm]
u
=
-
\nabla p
-
(0,\rho),
\qquad
\nabla
\cdot
u
=
0.
\end{cases}
\end{equation}
Here, $\rho$ is the density, $u=(u^1,u^2)$ is the velocity field, and $p$ is the pressure.

In this paper, we study the asymptotic stability of the steady state
\begin{align}\label{eq:steady}
e_s(x)
=
-
x_2,
\qquad
u_s
=
0,
\qquad
p_s(x)
=
\frac{x_2^2}{2},
\qquad
x
\in
\mathbb R^2.
\end{align}
To this end, we write $\rho=e_s+\theta$ and study the evolution
of the perturbation $\theta$.
A direct calculation shows that $\theta$ satisfies
\begin{equation}
\label{equ1.2}
\begin{cases}
\partial_t \theta
+
u
\cdot
\nabla \theta
=
\mathcal R_1^2
\theta,
\\[2mm]
u
=
\bigl(
-
\mathcal R_1
\mathcal R_2
\theta,\,
\mathcal R_1^2
\theta
\bigr),
\end{cases}
\end{equation}
where $\mathcal R_j=\partial_j(-\Delta)^{-1/2}$,
$j=1,2$, denotes the Riesz transforms.

\subsection{Known Results}

C\'ordoba, Gancedo, and Orive~\cite{CGO2007} established local
existence and uniqueness for regular solutions of the IPM
equation~\eqref{equ1.1}, together with continuation criteria.
Elgindi~\cite{Elgindi2017} first established the global stability of
the steady state $e_s(x)=-x_2$ in $\mathbb R^2$ and $\mathbb T^2$
for sufficiently small perturbations in $H^m\cap W^{4,1}$ and
$H^m$, respectively, with integers $m\ge20$.
In $\mathbb R^2$, he obtained algebraic decay of the density
perturbation and velocity in lower-order Sobolev norms, while on
$\mathbb T^2$, he established decay of the nonzero horizontal
Fourier modes of the density perturbation, as well as the velocity.
Castro, C\'ordoba, and Lear~\cite{CCL2019} extended the global
stability result to the periodic channel $\mathbb T\times[-1,1]$
for sufficiently small perturbations in $H^m$, with integers
$m\ge10$, subject to suitable boundary compatibility conditions.
Bianchini, Crin-Barat, and Paicu~\cite{BCP2024} established global
asymptotic stability in $\mathbb R^2$ for sufficiently small
perturbations in $\dot H^{1-\tau}\cap\dot H^s$, where
$0<\tau<1$ and $s\ge3+\tau$.
In particular, their result applies to small $H^s$ perturbations
for every $s>3$.
Jo and Kim~\cite{JK2024} independently established asymptotic
stability in the same $H^s$ regularity range.
Park~\cite{Park2024} established global stability in the periodic
channel $\mathbb T\times(0,1)$ for sufficiently small perturbations
in $H^m$, with integers $m\ge3$, subject to suitable boundary
compatibility conditions, by exploiting the potential energy
structure.
Building on Park's potential-energy approach, Bianchini, Jo, Park,
and Wang~\cite{BJPW2025} established global asymptotic stability
of $e_s(x)=-x_2$ (and, more generally, sufficiently smooth,
uniformly decreasing stratified states) on
$\mathbb T\times\mathbb R$ for sufficiently small $H^s$
perturbations, with every real $s>2$.
A key ingredient of~\cite{BJPW2025} is the establishment of
sufficiently strong time-decay estimates for the solutions.
In $\mathbb R^2$, establishing decay of the solutions under the
sole assumption of small $H^s$ initial data, without additional
$L^1$ or negative-order Sobolev assumptions, is a delicate issue.
Even when such decay can be established, it is generally slower
than on $\mathbb T\times\mathbb R$.
Therefore, extending the approach of~\cite{BJPW2025} to
$\mathbb R^2$ presents substantial difficulties.
Consequently, global asymptotic stability in $\mathbb R^2$ for
small $H^s$ perturbations in the range $2<s\le3$ remained open.

On the other hand, instability phenomena have also been observed
for the steady state $e_s(x)=-x_2$ at lower Sobolev regularities.
Kiselev and Yao~\cite{KY2023} established nonlinear instability
of this steady state (and, more generally, of arbitrary smooth
stratified steady states) in the periodic channel
$\mathbb T\times[-\pi,\pi]$.
For every $\gamma>0$, they constructed smooth initial perturbations
arbitrarily small in $H^{2-\gamma}$ such that the corresponding
density perturbations exhibit unbounded growth in $H^s$ norms
for every $s>1$, provided the solutions remain globally smooth.
In the endpoint space $H^2$, which is also a critical space for
the IPM equation, Bianchini, C\'ordoba, and
Mart\'inez-Zoroa~\cite{BCM2025} established instability of the
steady state $e_s(x)=-x_2$ in $\mathbb R^2$.
They constructed arbitrarily small $H^2$ initial perturbations
for which the corresponding density perturbations instantaneously
lose $H^2$ regularity.
For other stratified backgrounds $\rho_s(x_2)$ that increase
with $x_2$, Banihashemi, Mohammadkhani, and
Nguyen~\cite{BMN2026} established some nonlinear instability results.

These results naturally raise the following questions in
$\mathbb R^2$: whether the steady state $e_s(x)=-x_2$ is globally asymptotically stable under sufficiently small $H^s$ perturbations for $2<s\le3$. 
\emph{Remarkably, it is not clear whether the corresponding density perturbations will
converge to zero in $H^s$}. Regarding instability, the result of Kiselev and Yao~\cite{KY2023}
for perturbations of this steady state in a periodic channel does not describe the growth of $H^s$ norms for $0<s\le1$.
Furthermore, the local-in-time instability in
$H^s(\mathbb R^2)$ for $0<s<2$ is still unknown. The present paper focuses on the stability problem.
In a forthcoming paper \cite{LL2026}, we will establish local-in-time instability
by proving norm inflation of the density perturbation in
$H^s(\mathbb R^2)$ for all $0<s<2$.

\subsection{Main theorem}

Our main result reads as follows.

\begin{Theorem}\label{thm:main}
For any $s>2$, there exist two positive constants $\ep_s$
and $C_s$, depending only on $s$, such that for any initial datum
$\theta_0$ satisfying $\|\theta_0\|_{H^s}\le\ep_s$,
the system \eqref{equ1.2} admits a unique global solution
$\theta\in C([0,\infty);H^s)$.
Moreover, the solution satisfies
\begin{equation}
\label{eq:global-bound}
\sup_{t\ge0}
\|\theta(t)\|_{H^s}^{2}
+
\int_0^\infty
\|\mathcal R_1\theta(t)\|_{H^s}^{2}\,
dt
\le
C_s
\|\theta_0\|_{H^s}^{2},
\end{equation}
and
\begin{equation}
\label{eq:main-convergence}
\lim_{t\to\infty}
\|\theta(t)\|_{H^s}
=
0.
\end{equation}
\end{Theorem}

%\begin{Remark}
%The global stability arguments
%in~\cite{BCP2024,BJPW2025,CCL2019,Elgindi2017,JK2024,Park2024}
%involve control of $\|\p_2 u_2\|_{L^1_tL^\infty_x}$.
%Our weighted energy estimate closes directly for
%small $H^s$ initial data, with any $s>2$, without requiring
%this control.
%Moreover, the convergence results for the density in these works
%are limited to lower-order Sobolev norms.
%We establish strong convergence of the density perturbation to zero
%in the highest-order $H^s(\mathbb R^2)$ norm,
%at the same regularity as the initial datum.
%\end{Remark}

%\begin{Remark}
%Almost simultaneously, Bianchini et al. \cite{BCSW2026} independently established
%global stability for sufficiently small initial perturbations
%$\theta_0$ in $\dot B^1_{2,1}\cap\dot B^2_{2,1}$, as well as %in
%$H^s$ with $s>2$. However, they did not establish the convergence stated in \eqref{eq:main-convergence}. Their proof relies on control of
%$\|\p_2 u_2\|_{L^1_tL^\infty_x}$.
%In comparison, our weighted energy argument avoids this control
%and yields a simpler proof of global stability.
%\end{Remark}

After the completion of this manuscript, we noticed on arXiv an independent work by Bianchini et al. \cite{BCSW2026}, which also studied the asymptotic stability of the IPM equation on $\mathbb{R}^2$ around the steady state \eqref{eq:steady}. However, there are several distinct features between \cite{BCSW2026} and our main Theorem \ref{thm:main}. More precisely, 
\begin{itemize}
    \item{
    As in many previous studies \cite{BCP2024,BJPW2025,CCL2019,Elgindi2017,JK2024,Park2024}, the approach in Bianchini et al. \cite{BCSW2026} requires a uniform bound on $\|\Grad u^2\|_{L^1_t L^\infty_x}$, which is obtained through a detailed and technically involved analysis in anisotropic homogeneous Besov spaces.
By exploiting the structure of the equation \eqref{equ1.2},
we construct a weighted energy that eliminates the need
for this control and yields a simplified proof of global stability.
    }
    %\item{
   % Under the additional assumption that $\theta_0 \in W^{4,1}(\mathbb{R}^2)$, Bianchini et al. \cite{BCSW2026} proved $H^2$ decay estimates of global solutions.
    %in $H^2(\mathbb{R}^2)$. 
   % Together with the uniform $H^s$ bound, these estimates
%imply that $\theta$ converges to zero in $H^{s'}(\mathbb R^2)$
%for every $0\le s'<s$.
    %in which the equation~\eqref{equ1.2} is ill-posed. 
    %In contrast, our decay estimate \eqref{eq:main-convergence} holds for the highest-order Sobolev norm $H^s(\mathbb{R}^2)$, 
%without any additional assumptions beyond smallness
%of the initial data in $H^s$.
    %in which the equation~\eqref{equ1.2} is well-posed. 
   % }
    \item{
    It is easily shown %see Lemma \ref{lem:l2-convergence}, 
    that if
$z\in C([0,\infty);L^2(\mathbb R^2))$ solves
\[
\partial_t z=\mathcal R_1^2z+f,
\]
with $f\in L^1(0,\infty;L^2(\mathbb R^2))$, then
\[
\lim_{t\to\infty}\|z(t)\|_{L^2}=0.
\label{eq:remark-l1-forcing-convergence}
\]
%This observation will be used to prove \eqref{eq:main-convergence}. 
Both the global stability results in \cite{BCP2024,BJPW2025,CCL2019,Elgindi2017,JK2024,Park2024}
and those of Bianchini et al. \cite{BCSW2026}
involve the estimates for
$\| u^2\|_{L^1(0,\infty;L^\infty)}$, which implies $u\cdot \Grad \theta \in  L^1(0,\infty;L^2(\mathbb R^2))$.
Therefore, these results also imply the strong convergence of $\theta$ to zero  in the $L^2$ norm, and  by interpolation, in lower-order
$H^{s'}$ norms for all $0\leq s'<s$.
However, the convergence in the highest-order $H^s$ norm
does not follow directly.
In the present work, we establish
the highest-order convergence stated in \eqref{eq:main-convergence}
for the first time. Remarkably, the control of $\|u^2\|_{L^1(0,\infty;L^\infty)}$ is not needed.
    }
\end{itemize}

\subsection{Difficulties and strategy of the proof}

The standard approach to prove \eqref{eq:global-bound}
consists of two steps.
The first step is to establish an energy estimate of the form
\begin{equation*}
\begin{aligned}
&\|\theta\|_{L^\infty(0,t;H^s)}^2
+\|\mathcal R_1\theta\|_{L^2(0,t;H^s)}^2\\
&\quad\leq C\Bigl(
\|\theta_0\|_{H^s}^2
+\|\nabla u^2\|_{L^1(0,t;L^\infty)}
 \|\theta\|_{L^\infty(0,t;H^s)}^2\\
&\qquad\qquad+
\bigl[
\|\theta\|_{L^\infty(0,t;H^s)}
+\|\mathcal R_1\theta\|_{L^2(0,t;H^s)}
\bigr]^3
\Bigr).
\end{aligned}
\end{equation*}
The second step is to control
$\|\nabla u^2\|_{L^1(0,t;L^\infty)}$.
There are two typical ways to obtain this control:
deriving time-decay estimates for $u$,
or estimating $\|u\|_{L^1(0,t;B)}$
in a suitable Besov or anisotropic Besov space $B$.
For the spatial domain and initial-data assumptions
considered here, the latter approach is more natural,
but it entails greater technical complexity.

Our idea is to eliminate the need to control
$\|\nabla u^2\|_{L^1(0,t;L^\infty)}$.
To this end, we exploit the structure of equation~\eqref{equ1.2}
to construct a weighted energy $D_s(t)$,
which is equivalent to $E_s(t)$ provided that
$\|\partial_2\theta\|_{L^\infty(0,t;L^\infty)}<1$.
This construction is presented in
Section~\ref{sec:weighted-construction}.

The main difficulty in estimating $D_s$ remains
the treatment of anisotropic commutators.
Rather than estimating these commutators directly,
we study the integral terms containing them.
We express these integrals in Fourier variables
and carry out refined computations in frequency space,
obtaining suitable bounds for $D_s$ throughout
the regularity range $s>2$;
see Section~\ref{sec:weighted-energy-estimates}.
These estimation methods may also be useful in the analysis of related anisotropic problems.

For the convergence in \eqref{eq:main-convergence},
the main difficulty is that the energy bounds encoded
by $D_s$ do not provide any $L^1$-in-time dissipation bound.
We overcome this difficulty in three steps.
First, we establish an $L^2$-convergence lemma
that does not require such a bound;
see Lemma~\ref{lem:l2-convergence}.
Second, we apply a density-coordinate transformation
to equation~\eqref{equ1.2},
which allows us to use this lemma to establish
$L^2$ convergence;
see Proposition~\ref{prop:density-convergence}.
Third, using Proposition~\ref{Pro1},
we obtain uniform estimates for $D_{s'}$,
with constants independent of $s'$
as $s'$ approaches $s$.
A limiting argument then yields convergence
in the highest-order $H^s$ norm,
as stated in \eqref{eq:main-convergence}.

\section{The weighted energy \texorpdfstring{$D_s$}{Ds}}

\subsection{Construction of \texorpdfstring{$D_s$}{Ds}}\label{sec:weighted-construction}

This subsection is devoted to the construction of a weighted energy
that avoids the need to control
$\|\nabla u^2\|_{L^1_tL^\infty_x}$.
When deriving $H^s$ energy estimates for
equation~\eqref{equ1.2}, we encounter the term
\[
\int_{\mathbb R^2}
[\Lambda^s,u^2]\partial_2\theta\,\Lambda^s\theta\,dx.
\]
Here, for any $f$, we define
\[
\widehat f(\xi)
=
\frac{1}{2\pi}
\int_{\mathbb R^2}e^{-ix\cdot\xi}f(x)\,dx,
\qquad
\Lambda^sf
=
\mathcal F^{-1}\bigl(|\xi|^s\widehat f(\xi)\bigr),
\]
where $\mathcal F^{-1}$ denotes the inverse Fourier transform.
Applying a standard commutator estimate, we obtain
\[
\left|
\int_{\mathbb R^2}
[\Lambda^s,u^2]\partial_2\theta\,\Lambda^s\theta\,dx
\right|
\leq
C\left[
\|\nabla u^2\|_{L^\infty}\|\theta\|_{H^s}^2
+
\|\mathcal R_1\theta\|_{H^s}\|\theta\|_{H^s}^2
\right].
\]
Through more refined computations, the second term can be improved to
$
\|\mathcal R_1\theta\|_{H^s}^2\|\theta\|_{H^s}.
$
However, the first term $\|\nabla u^2\|_{L^\infty}\|\theta\|_{H^s}^2$ remains.

We therefore use the norm equivalence
\[
\|\theta\|_{H^s}
\approx
\|\theta\|_{L^2}
+\|\Lambda_1^s\theta\|_{L^2}
+\|\Lambda_2^s\theta\|_{L^2}
\]
and estimate the three terms separately.
Here, for any $f$, we define
\[
\Lambda_i^sf
:=
\mathcal F^{-1}\bigl(|\xi^{(i)}|^s\widehat f(\xi)\bigr),
\qquad i=1,2,
\]
where $\xi=(\xi^{(1)},\xi^{(2)})$.
The difficulty arises when estimating $\Lambda_2^s\theta$. However,
a direct computation reveals that the problematic term involves
only $\partial_2u^2$, rather than the full gradient $\nabla u^2$.
Specifically, applying $\Lambda_2^s$ to \eqref{equ1.2}, multiplying by $\Lambda_2^s\theta$, and integrating over $\mathbb{R}^2$, we obtain
\[
\begin{aligned}
\frac12
\frac{d}{dt}
\norm{\Lambda_2^s\theta}_{L^2}^2
+
\norm{\mathcal R_1\Lambda_2^s\theta}_{L^2}^2
&
=
-
\int_{\mathbb{R}^2}
[\Lambda_2^s,u^1]
\partial_1 \theta\,
\Lambda_2^s \theta
\dd x
\\
&\quad
-
\int_{\mathbb{R}^2}
[\Lambda_2^s,u^2]
\partial_2 \theta\,
\Lambda_2^s \theta
\dd x
\\
&
:=
-
P_1
-
P_2.
\end{aligned}
\]
%Here, for any $f$, the directional derivatives are defined by
%\[
%\Lambda_j^s f
%=
%\mathcal F^{-1}
%\bigl(
%|\xi^{(j)}|^s
%\widehat f(\xi)
%\bigr),
%\qquad
%j
%=
%1,
%2,
%\]
%where $\xi=(\xi^{(1)},\xi^{(2)})$ and $\mathcal F^{-1}$ denotes the inverse Fourier transform.

Writing $P_2$ in Fourier variables and using the first-order expansion of $|\cdot|^s$, we obtain
\[
\begin{aligned}
P_2
={}
&
\frac{i}{2\pi}
\int_{\mathbb{R}^2}
\int_{\mathbb{R}^2}
\Bigl[
|\xi_2^{(2)}|^s
-
|\xi_1^{(2)}|^s
-
s
\bigl(
\xi_2^{(2)}
-
\xi_1^{(2)}
\bigr)
\xi_1^{(2)}
|\xi_1^{(2)}|^{s-2}
\Bigr]
\\
&\qquad
\times
\widehat{u^2}(\xi_2-\xi_1)\,
\xi_1^{(2)}
\widehat\theta(\xi_1)\,
|\xi_2^{(2)}|^s
\widehat\theta(-\xi_2)
\dd\xi_1
\dd\xi_2
\\
&
+
\frac{s}{2\pi}
\int_{\mathbb{R}^2}
\int_{\mathbb{R}^2}
\widehat{u^2}(\xi_2-\xi_1)\,
i
\bigl(
\xi_2^{(2)}
-
\xi_1^{(2)}
\bigr)\,
|\xi_1^{(2)}|^s
\widehat\theta(\xi_1)\,
|\xi_2^{(2)}|^s
\widehat\theta(-\xi_2)
\dd\xi_1
\dd\xi_2
\\[4pt]
={}
&
\frac{i}{2\pi}
\int_{\mathbb{R}^2}
\int_{\mathbb{R}^2}
\Bigl[
|\xi_2^{(2)}|^s
-
|\xi_1^{(2)}|^s
-
s
\bigl(
\xi_2^{(2)}
-
\xi_1^{(2)}
\bigr)
\xi_1^{(2)}
|\xi_1^{(2)}|^{s-2}
\Bigr]
\\
&\qquad
\times
\widehat{u^2}(\xi_2-\xi_1)\,
\xi_1^{(2)}
\widehat\theta(\xi_1)\,
|\xi_2^{(2)}|^s
\widehat\theta(-\xi_2)
\dd\xi_1
\dd\xi_2
\\
&
+
s
\int_{\mathbb{R}^2}
\partial_2 u^2\,
(\Lambda_2^s\theta)^2
\dd x.
\end{aligned}
\]
To handle the last term, we consider the equation
\begin{equation}
\label{eq:q}
\partial_t \Lambda_2^s \theta
+
u
\cdot
\nabla \Lambda_2^s \theta
=
\mathcal R_1^2
\Lambda_2^s \theta
-
s
\partial_2 u^2
\Lambda_2^s \theta
-
R_h
-
R_v,
\end{equation}
where
\[
R_h
=
[\Lambda_2^s,u^1]
\partial_1 \theta,
\qquad
R_v
=
[\Lambda_2^s,u^2]
\partial_2 \theta
-
s
\partial_2 u^2
\Lambda_2^s \theta.
\]

To obtain an equation involving $\partial_2u^2$, we differentiate \eqref{equ1.2} in the $x_2$ direction. This gives
\[
\begin{aligned}
\partial_t \partial_2 \theta
+
u
\cdot
\nabla \partial_2 \theta
&
=
\partial_2 u^2
-
\partial_2 u^2
\partial_2 \theta
-
\partial_2 u^1
\partial_1 \theta
\\
&
=
\partial_2 u^2(1-\partial_2\theta)
-
\partial_2 u^1
\partial_1 \theta.
\end{aligned}
\]
The factor $1-\partial_2\theta$ suggests considering
\[
\partial_t
(1-\partial_2\theta)
+
u\cdot\nabla
(1-\partial_2\theta)
=
-
(1-\partial_2\theta)
\partial_2 u^2
+
\partial_2 u^1
\partial_1 \theta.
\]
Assuming that $\norm{\partial_2\theta}_{L^\infty}<1$, we introduce
\[
\theta_s
:=
(1-\partial_2\theta)^{-2s}
\]
to match the coefficient $s$ in \eqref{eq:q}. It follows that
\begin{equation}
\label{eq:weight-transport}
\partial_t \theta_s
+
u
\cdot
\nabla \theta_s
=
2s
\theta_s
\partial_2 u^2
-
2s
\frac{\theta_s}{1-\partial_2\theta}
\partial_2 u^1
\partial_1 \theta.
\end{equation}

We multiply \eqref{eq:q} by $\Lambda_2^s\theta\,\theta_s$ and \eqref{eq:weight-transport} by $(\Lambda_2^s\theta)^2/2$, integrate over $\mathbb{R}^2$, and add the resulting identities. The terms containing $s\partial_2u^2(\Lambda_2^s\theta)^2\theta_s$ cancel. Hence,
\begin{equation}
\label{eq:weighted-vertical}
\begin{aligned}
&
\frac12
\frac{d}{dt}
\int_{\mathbb{R}^2}(\Lambda_2^s\theta)^2
\theta_s
\dd x
+
\norm{\mathcal R_1\Lambda_2^s\theta}_{L^2}^2
\\
&\quad
=
\int_{\mathbb{R}^2}
\mathcal R_1^2
\Lambda_2^s \theta\,
\Lambda_2^s \theta(\theta_s-1)
\dd x
%\\
%&\qquad
-
s
\int_{\mathbb{R}^2}
\frac{\theta_s}{1-\partial_2\theta}
\partial_2 u^1
\partial_1 \theta(\Lambda_2^s\theta)^2
\dd x
\\
&\qquad
-
\int_{\mathbb{R}^2}
R_h
\Lambda_2^s \theta\,
\theta_s
\dd x
-
\int_{\mathbb{R}^2}
R_v
\Lambda_2^s \theta\,
\theta_s
\dd x.
\end{aligned}
\end{equation}

Next, multiplying \eqref{equ1.2} by $\theta$ and integrating, and then applying $\Lambda_1^s$ to \eqref{equ1.2}, multiplying by $\Lambda_1^s\theta$, and integrating, we combine the two resulting identities with \eqref{eq:weighted-vertical} to obtain
\begin{equation}
\label{eq:full-energy}
\begin{aligned}
&
\frac12
\frac{d}{dt}
\left(
\norm{\theta}_{L^2}^2
+
\norm{\Lambda_1^s\theta}_{L^2}^2
+
\int_{\mathbb{R}^2}(\Lambda_2^s\theta)^2
\theta_s
\dd x
\right)
\\
&\quad
+
\norm{\mathcal R_1\theta}_{L^2}^2
+
\norm{\mathcal R_1\Lambda_1^s\theta}_{L^2}^2
+
\norm{\mathcal R_1\Lambda_2^s\theta}_{L^2}^2
\\
&
=
\int_{\mathbb{R}^2}
\mathcal R_1^2
\Lambda_2^s \theta\,
\Lambda_2^s \theta(\theta_s-1)
\dd x
%\\
%&\quad
-
\int_{\mathbb{R}^2}
R_h
\Lambda_2^s \theta\,
\theta_s
\dd x
-
\int_{\mathbb{R}^2}
R_v
\Lambda_2^s \theta\,
\theta_s
\dd x
\\
&\quad
-
s
\int_{\mathbb{R}^2}
\frac{\theta_s}{1-\partial_2\theta}
\partial_2 u^1
\partial_1 \theta(\Lambda_2^s\theta)^2
\dd x
\\
&\quad
-
\int_{\mathbb{R}^2}
[\Lambda_1^s,u^1]
\partial_1 \theta\,
\Lambda_1^s \theta
\dd x
-
\int_{\mathbb{R}^2}
[\Lambda_1^s,u^2]
\partial_2 \theta\,
\Lambda_1^s \theta
\dd x.
\end{aligned}
\end{equation}

\begin{samepage}
For $s>2$, $0\le T<\infty$, and $0\le T_1<T_2\le\infty$, we define
\[
\begin{aligned}
P_s(T)
:=
\Big(
&
\|\theta(T)\|_{L^2}^2
+
\|\Lambda_1^s\theta(T)\|_{L^2}^2
+
\|\sqrt{\theta_s(T)}\,\Lambda_2^s\theta(T)\|_{L^2}^2
\Big)^{1/2},
\end{aligned}
\]
\[
\begin{aligned}
D_s(T_1,T_2)
:=
\Big(
&
\|P_s\|_{L^\infty(T_1,T_2)}^2
+
\|\mathcal R_1\theta\|_{L^2(T_1,T_2;L^2)}^2
+
\|\mathcal R_1\Lambda_1^s\theta\|_{L^2(T_1,T_2;L^2)}^2
\\
&
+
\|\mathcal R_1\Lambda_2^s\theta\|_{L^2(T_1,T_2;L^2)}^2
\Big)^{1/2},
\end{aligned}
\]
\[
\begin{aligned}
E_s(T_1,T_2)
:={}
&
\|\theta\|_{L^\infty(T_1,T_2;H^s)}
+
\|\mathcal R_1\theta\|_{L^2(T_1,T_2;H^s)}.
\end{aligned}
\]
\end{samepage}
For simplicity, we write $D_s(0,T)$ and $E_s(0,T)$ as $D_s(T)$ and $E_s(T)$, respectively.
If
\[
\|\partial_2\theta(T)\|_{L^\infty}
\le
\alpha
<
1,
\qquad
\|\partial_2\theta\|_{L^\infty(T_1,T_2;L^\infty)}
\le
\alpha
<
1,
\]
then
\begin{equation}
\label{eq:Ps-equivalence}
C_{\alpha,s}^{-1}
\|\theta(T)\|_{H^s}
\le
P_s(T)
\le
C_{\alpha,s}
\|\theta(T)\|_{H^s},
\end{equation}
\begin{equation}
\label{eq:Ds-Es-equivalence}
C_{\alpha,s}^{-1}
E_s(T_1,T_2)
\le
D_s(T_1,T_2)
\le
C_{\alpha,s}
E_s(T_1,T_2),
\end{equation}
where
\begin{equation}
C_{\alpha,s}
=
\max
\left\{
(1-\alpha)^{-s},
3^{\frac{s-1}{2}}
\sqrt{1+(1+\alpha)^{2s}}
\right\}.
\label{eq:energy-equivalence-constant}
\end{equation}

\subsection{Estimates for \texorpdfstring{$D_s$}{Ds}}
\label{sec:weighted-energy-estimates}

This subsection is devoted to the following proposition.
\begin{Proposition}\label{Pro1}
Let $s>2$, $0<\alpha<1$, and $0\le T_1<T_2<\infty$.
Suppose that $\theta$ is a regular solution of \eqref{equ1.2} on $[T_1,T_2]$ satisfying
\[
\|\theta\|_{L^\infty(T_1,T_2;H^s)}
\leq
2
\alpha
\sqrt{\pi(s-2)},
\qquad
\|\mathcal R_1\theta\|_{L^2(T_1,T_2;H^s)}
<
\infty.
\]
Then the following estimates hold:
\[
\begin{aligned}
D_s(T_1,T_2)^2
\leq{}
&
\frac32
P_s(T_1)^2
+
\Big[
6s(1-\alpha)^{-(2s+1)}
\alpha^2
+
3s^2
2^{s+2}
C^*_{\alpha, s}
\Big]
\times
\|\mathcal R_1\theta\|_{L^2(T_1,T_2;H^s)}^2,
\end{aligned}
\]
\[
\begin{aligned}
P_s(T_1)^2
\leq{}
&
P_s(T_2)^2
+
\Big[
2
+
4s(1-\alpha)^{-(2s+1)}
\alpha^2
+
s^2
2^{s+3}
C^*_{\alpha, s}
\Big]
\times
\|\mathcal R_1\theta\|_{L^2(T_1,T_2;H^s)}^2.
\end{aligned}
\]
where
\begin{align}\label{new1}
C^*_{\alpha, s}
:=
\sum_{n=1}^{\infty}
e^{2s-1}
[
n^{2s-1}
+
(n+1)^{2s-1}
]
(n+1)^2
\alpha^n.
\end{align}
\end{Proposition}

\subsubsection{Elementary estimates and Fourier notation}
In this subsection, we introduce some basic estimates that will be used later. For $s>2$ and $a,b\in\mathbb{R}$,
\begin{align}
\int_{\mathbb{R}^2}
\frac{1}{\langle x\rangle^s}
\dd x
&
=
\frac{2\pi}{s-2},
\qquad
\langle
x
\rangle
=
(1+|x|^2)^{1/2},
\label{eq:elementary-integral}
\\
\bigl|
|a+b|^s
-
|a|^s
\bigr|
&
\le
s2^{s-1}(|b||a|^{s-1}+|b|^s),
\label{eq:first-diff}
\\
\bigl|
|a+b|^s
-
|a|^s
-
sa
|a|^{s-2}
b
\bigr|
&
\le
s(s-1)
2^{s-2}(|b|^2|a|^{s-2}+|b|^s).
\label{eq:second-diff}
\end{align}
For $N\ge1$, we have
\[
\begin{aligned}
&
\int_{\mathbb{R}^2}
f_1(x)
\cdots
f_N(x)
f_{N+1}(x)
\dd x
\\
&\quad
=
(2\pi)^{-(N-1)}
\int_{(\mathbb{R}^2)^N}
\widehat f_1(\xi_1)
\widehat f_2(\xi_2-\xi_1)
\cdots
\widehat f_N(\xi_N-\xi_{N-1})
\widehat f_{N+1}(-\xi_N)
\dd\xi_1
\cdots
\dd\xi_N.
\end{aligned}
\]
For convenience, we write
\[
\begin{aligned}
&
\int_{(\mathbb{R}^2)^N}
\widehat f_1(\xi_1)
\widehat f_2(\xi_2-\xi_1)
\cdots
\widehat f_N(\xi_N-\xi_{N-1})
\widehat f_{N+1}(-\xi_N)\,
d
\xi_1
\cdots
d
\xi_N
\end{aligned}
\]
as
\[
\int_{\xi_1+\cdots+\xi_{N+1}=0}
\widehat f_1(\xi_1)
\cdots
\widehat f_{N+1}(\xi_{N+1})\,
d
\xi_1
\cdots
d
\xi_N.
\]
Young's convolution inequality and H\"older's inequality yield
\begin{equation}
\label{eq:young}
\left|
\int_{\xi_1+\cdots+\xi_{N+1}=0}
\prod_{j=1}^{N+1}
\widehat f_j(\xi_j)
\dd\xi_1
\cdots
\dd\xi_N
\right|
\le
\left(
\prod_{k=1}^{N-1}
\norm{\widehat f_{i_k}}_{L^1}
\right)
\norm{\widehat f_{i_N}}_{L^p}
\norm{\widehat f_{i_{N+1}}}_{L^q},
\end{equation}
where $1\le p,q\le\infty$, $1/p+1/q=1$, and $(i_1,\ldots,i_{N+1})$ is a permutation of $(1,\ldots,N+1)$.

From \eqref{eq:elementary-integral}, we obtain
\begin{align}
\int_{\mathbb{R}^2}
|\xi|\,
|\widehat\theta(\xi)|
\dd\xi
&
\le
\sqrt{\frac{\pi}{s-2}}
\norm{\theta}_{H^s},
\label{eq:L1theta}
\\
\int_{\mathbb{R}^2}
|\xi^{(1)}|\,
|\widehat\theta(\xi)|
\dd\xi
&
\le
\sqrt{\frac{\pi}{s-2}}
\norm{\mathcal R_1\theta}_{H^s},
\label{eq:L1horizontal}
\\
\left(
\int_{\mathbb{R}^2}
|\xi^{(1)}|^2
|\xi|^{2(s-1)}
|\widehat\theta(\xi)|^2
\dd\xi
\right)^{1/2}
&
\le
\norm{\mathcal R_1\theta}_{H^s}.
\label{eq:L2horizontal}
\end{align}
For $|a|<1$, we use the expansion
\[
(1-a)^{-2s}
=
\sum_{n=0}^{\infty}
c_n
a^n,
\quad
c_0
=
1,
\quad
c_n
=
\frac{(2s)(2s+1)\cdots(2s+n-1)}{n!}
=
\prod_{k=1}^n
\left(
1
+
\frac{2s-1}{k}
\right).
\]
A direct calculation gives
\[
\log
c_n
\le
(2s-1)
\sum_{k=1}^n
\frac1k
\le
(2s-1)(1+\log n),
\qquad
n
\ge
1.
\]
Consequently,
\begin{equation}
\label{eq:cn-bound}
|c_n|
\le
e^{2s-1}(n+1)^{2s-1},
\qquad
n
\ge
0.
\end{equation}

\subsubsection{The proof of Proposition~\ref{Pro1}}

\begin{proof}
From \eqref{eq:L1theta}, we get
\[
\begin{aligned}
\norm{\partial_2\theta}_{L^\infty(T_1,T_2;L^\infty)}
&
\le
\frac{1}{2\pi}
\norm{\widehat{\partial_2\theta}}_{L^\infty(T_1,T_2;L^1)}
\le
\frac{\norm{\theta}_{L^\infty(T_1,T_2;H^s)}}{2\sqrt{\pi(s-2)}}
\le
\alpha.
\end{aligned}
\]
Using \eqref{eq:full-energy}, we obtain
\begin{equation}
\label{eq:Ds-source-bound}
D_s(T_1,T_2)^2
\leq
\frac32
P_s(T_1)^2
+
3
\sum_{j=1}^{6}
I_j(T_1,T_2).
\end{equation}
\begin{equation}
\label{eq:Ps-reverse-source-bound}
\begin{aligned}
P_s(T_1)^2
\leq{}
&
P_s(T_2)^2
+
2
\|\mathcal R_1\theta\|_{L^2(T_1,T_2;L^2)}^2
+
2
\|\mathcal R_1\Lambda_1^s\theta\|_{L^2(T_1,T_2;L^2)}^2
\\
&
+
2
\|\mathcal R_1\Lambda_2^s\theta\|_{L^2(T_1,T_2;L^2)}^2
+
2
\sum_{j=1}^{6}
I_j(T_1,T_2).
\end{aligned}
\end{equation}
where
\[
\begin{aligned}
I_1(T_1,T_2)
&
:=
s
\int_{T_1}^{T_2}
\left|
\int_{\mathbb{R}^2}
\frac{\theta_s}{1-\partial_2\theta}
\partial_2 u^1
\partial_1 \theta
 (\Lambda_2^s\theta)^2
\dd x
\right|
\dd\tau,
\\
I_2(T_1,T_2)
&
:=
\int_{T_1}^{T_2}
\left|
\int_{\mathbb{R}^2}
[\Lambda_1^s,u^1]
\partial_1 \theta\,
\Lambda_1^s \theta
\dd x
\right|
\dd\tau,
\\
I_3(T_1,T_2)
&
:=
\int_{T_1}^{T_2}
\left|
\int_{\mathbb{R}^2}
[\Lambda_1^s,u^2]
\partial_2 \theta\,
\Lambda_1^s \theta
\dd x
\right|
\dd\tau,
\\
I_4(T_1,T_2)
&
:=
\int_{T_1}^{T_2}
\left|
\int_{\mathbb{R}^2}
\mathcal R_1^2
\Lambda_2^s \theta\,
\Lambda_2^s \theta(\theta_s-1)
\dd x
\right|
\dd\tau,
\\
I_5(T_1,T_2)
&
:=
\int_{T_1}^{T_2}
\left|
\int_{\mathbb{R}^2}
R_v
\Lambda_2^s \theta\,
\theta_s
\dd x
\right|
\dd\tau,
\\
I_6(T_1,T_2)
&
:=
\int_{T_1}^{T_2}
\left|
\int_{\mathbb{R}^2}
R_h
\Lambda_2^s \theta\,
\theta_s
\dd x
\right|
\dd\tau.
\end{aligned}
\]

\paragraph{Estimates for \texorpdfstring{$I_1$}{I1}.}
Using \eqref{eq:L1horizontal} and
$\norm{\partial_2\theta}_{L^\infty(T_1,T_2;L^\infty)}\le\alpha$, we obtain
\begin{equation}
\label{eq:I1}
\begin{aligned}
I_1(T_1,T_2)
&
\le
s
\left\|
\frac{\theta_s}{1-\partial_2\theta}
\right\|_{L^\infty(T_1,T_2;L^\infty)}
\norm{\partial_2u^1}_{L^2(T_1,T_2;L^\infty)}
\\[-2pt]
&\qquad
\times
\norm{\partial_1\theta}_{L^2(T_1,T_2;L^\infty)}
\norm{\Lambda_2^s\theta}_{L^\infty(T_1,T_2;L^2)}^2
\\
&
\le
\frac{s(1-\alpha)^{-(2s+1)}}{2\pi(s-2)}
\norm{\theta}_{L^\infty(T_1,T_2;H^s)}^2
\norm{\mathcal R_1\theta}_{L^2(T_1,T_2;H^s)}^2
\\
&
\le
2s(1-\alpha)^{-(2s+1)}
\alpha^2
\norm{\mathcal R_1\theta}_{L^2(T_1,T_2;H^s)}^2.
\end{aligned}
\end{equation}

\paragraph{Estimates for \texorpdfstring{$I_2$}{I2}.}
Writing the commutator in Fourier variables and using \eqref{eq:first-diff}, we obtain
\begin{align}
I_2(T_1,T_2)
&
\le
\frac{s2^{s-1}}{2\pi}
\int_{T_1}^{T_2}
\int_{\xi_1+\xi_2+\xi_3=0}
\bigl(
|\xi_1^{(1)}|
|\xi_2^{(1)}|^{s-1}
+
|\xi_1^{(1)}|^s
\bigr)
|\widehat{u^1}(\xi_1)|
\notag
\\*[-2pt]
&\qquad\qquad
\times
|\xi_2^{(1)}|
|\widehat\theta(\xi_2)|
|\xi_3^{(1)}|^s
|\widehat\theta(\xi_3)|
\dd\xi_1
\dd\xi_2
\dd\tau
\notag
\\
&
\le
\frac{2s2^{s-1}}{2\pi}
\int_{T_1}^{T_2}
\int_{\xi_1+\xi_2+\xi_3=0}
|\xi_1^{(1)}|
|\widehat\theta(\xi_1)|\,
|\xi_2^{(1)}|^s
|\widehat\theta(\xi_2)|
\notag
\\*[-2pt]
&\qquad\qquad
\times
|\xi_3^{(1)}|^s
|\widehat\theta(\xi_3)|
\dd\xi_1
\dd\xi_2
\dd\tau
\notag
\\
&
\le
\frac{s2^s}{2\sqrt{\pi(s-2)}}
\int_{T_1}^{T_2}
\norm{\theta(\tau)}_{H^s}
\norm{\mathcal R_1\theta(\tau)}_{H^s}^2
\dd\tau
\notag
\\
&
\le
s2^s
\frac{ \norm{\theta}_{L^\infty(T_1,T_2;H^s)}}{2\sqrt{\pi(s-2)}}
\norm{\mathcal R_1\theta}_{L^2(T_1,T_2;H^s)}^2
\notag
\\
&
\le
s2^s
\alpha
\norm{\mathcal R_1\theta}_{L^2(T_1,T_2;H^s)}^2.
\label{eq:I2}
\end{align}

\paragraph{Estimates for \texorpdfstring{$I_3$}{I3}.}
Similarly, we have
\begin{align}
I_3(T_1,T_2)
&
\le
\frac{s2^{s-1}}{2\pi}
\int_{T_1}^{T_2}
\int_{\xi_1+\xi_2+\xi_3=0}
\bigl(
|\xi_1^{(1)}|
|\xi_2^{(1)}|^{s-1}
+
|\xi_1^{(1)}|^s
\bigr)
|\widehat{u^2}(\xi_1)|
\notag
\\*[-2pt]
&\qquad\qquad
\times
|\xi_2^{(2)}|
|\widehat\theta(\xi_2)|
|\xi_3^{(1)}|^s
|\widehat\theta(\xi_3)|
\dd\xi_1
\dd\xi_2
\dd\tau
\notag
\\
&
\le
\frac{s2^{s-1}}{2\pi}
\int_{T_1}^{T_2}
\int_{\xi_1+\xi_2+\xi_3=0}
|\xi_1^{(1)}|
|\widehat\theta(\xi_1)|\,
|\xi_2|^s
|\widehat\theta(\xi_2)|
\notag
\\*[-2pt]
&\qquad\qquad
\times
|\xi_3^{(1)}|^s
|\widehat\theta(\xi_3)|
\dd\xi_1
\dd\xi_2
\dd\tau
\notag
\\
&\quad
+
\frac{s2^{s-1}}{2\pi}
\int_{T_1}^{T_2}
\int_{\xi_1+\xi_2+\xi_3=0}
|\xi_1^{(1)}|^s
|\widehat\theta(\xi_1)|\,
|\xi_2|
|\widehat\theta(\xi_2)|
\notag
\\*[-2pt]
&\qquad\qquad
\times
|\xi_3^{(1)}|^s
|\widehat\theta(\xi_3)|
\dd\xi_1
\dd\xi_2
\dd\tau
\notag
\\
&
\le
\frac{s2^s}{2\sqrt{\pi(s-2)}}
\int_{T_1}^{T_2}
\norm{\theta(\tau)}_{H^s}
\norm{\mathcal R_1\theta(\tau)}_{H^s}^2
\dd\tau
\notag
\\
&
\le
s2^s
\frac{\norm{\theta}_{L^\infty(T_1,T_2;H^s)}}{2\sqrt{\pi(s-2)}}
\norm{\mathcal R_1\theta}_{L^2(T_1,T_2;H^s)}^2
\notag
\\
&
\le
s2^s
\alpha
\norm{\mathcal R_1\theta}_{L^2(T_1,T_2;H^s)}^2.
\label{eq:I3}
\end{align}

\paragraph{Estimates for \texorpdfstring{$I_4$}{I4}.}
Expanding $\theta_s-1$, we obtain
\begin{align}
I_4(T_1,T_2)
&
\le
\sum_{n=1}^{\infty}
|c_n|
\int_{T_1}^{T_2}
\left|
\int_{\mathbb{R}^2}
\mathcal R_1^2
\Lambda_2^s \theta\,
\Lambda_2^s \theta
 (\partial_2\theta)^n
\dd x
\right|
\dd\tau
\notag
\\
&
\le
\sum_{n=1}^{\infty}
|c_n|(2\pi)^{-n}
\int_{T_1}^{T_2}
\int_{\xi_1+\cdots+\xi_{n+2}=0}
|\xi_1^{(1)}|^2
|\xi_1|^{s-2}
|\widehat\theta(\xi_1)|
|\xi_2|^s
|\widehat\theta(\xi_2)|
\notag
\\*[-2pt]
&\qquad\qquad
\times
\prod_{k=3}^{n+2}
|\xi_k|
|\widehat\theta(\xi_k)|
\dd\xi_1
\cdots
\dd\xi_{n+1}
\dd\tau.
\label{eq:I4-exp}
\end{align}
Using $|\xi_1^{(1)}|\le\sum_{j=2}^{n+2}|\xi_j^{(1)}|$ and symmetry among
$\xi_3,\ldots,\xi_{n+2}$, we obtain
\begin{align}
&
\int_{\xi_1+\cdots+\xi_{n+2}=0}
|\xi_1^{(1)}|^2
|\xi_1|^{s-2}
|\widehat\theta(\xi_1)|
|\xi_2|^s
|\widehat\theta(\xi_2)|
\prod_{k=3}^{n+2}
|\xi_k|
|\widehat\theta(\xi_k)|
\dd\xi_1
\cdots
\dd\xi_{n+1}
\notag
\\
&\quad
\le
\int_{\xi_1+\cdots+\xi_{n+2}=0}
|\xi_1^{(1)}|
|\xi_1|^{s-2}
|\widehat\theta(\xi_1)|
|\xi_2^{(1)}|
|\xi_2|^s
|\widehat\theta(\xi_2)|
\notag
\\*[-2pt]
&\qquad\qquad
\times
\prod_{k=3}^{n+2}
|\xi_k|
|\widehat\theta(\xi_k)|
\dd\xi_1
\cdots
\dd\xi_{n+1}
\notag
\\
&\qquad
+
n
\int_{\xi_1+\cdots+\xi_{n+2}=0}
|\xi_1^{(1)}|
|\xi_1|^{s-2}
|\widehat\theta(\xi_1)|
|\xi_2|^s
|\widehat\theta(\xi_2)|
\notag
\\*[-2pt]
&\qquad\qquad
\times
|\xi_3^{(1)}|
|\xi_3|
|\widehat\theta(\xi_3)|
\prod_{k=4}^{n+2}
|\xi_k|
|\widehat\theta(\xi_k)|
\dd\xi_1
\cdots
\dd\xi_{n+1}
\notag
\\
&\quad
=:
I_{4,1}^{(n)}
+
nI_{4,2}^{(n)}.
\label{eq:J12}
\end{align}
Let $|\xi_{i_1}|$ and $|\xi_{i_2}|$ be the largest and second-largest of
$|\xi_1|,\ldots,|\xi_{n+2}|$. Since $\xi_1+\cdots+\xi_{n+2}=0$,
\begin{equation}
\label{eq:two-largest}
|\xi_{i_1}|
\le
\sum_{j\ne i_1}
|\xi_j|
\le
(n+1)
|\xi_{i_2}|.
\end{equation}
We consider the following cases.

\textit{Case 1:} $i_1=2$, $i_2=1$. Using \eqref{eq:two-largest} and $|\xi_3|\le|\xi_1|$, we have
\[
\begin{aligned}
I_{4,1}^{(n)}
&
\le
(n+1)
\int_{\xi_1+\cdots+\xi_{n+2}=0}
|\xi_{1}^{(1)}|
|\xi_{1}|^{s-1}
|\widehat\theta(\xi_{1})|\,
|\xi_{2}^{(1)}|
|\xi_{2}|^{s-1}
|\widehat\theta(\xi_{2})|
\\[-2pt]
&\qquad
\times
\prod_{k=3}^{n+2}
|\xi_{k}|
|\widehat\theta(\xi_{k})|
\dd\xi_1
\cdots
\dd\xi_{n+1}.
\end{aligned}
\]
\[
\begin{aligned}
I_{4,2}^{(n)}
&
\le
\int_{\xi_1+\cdots+\xi_{n+2}=0}
|\xi_{1}^{(1)}|
|\xi_{1}|^{s-1}
|\widehat\theta(\xi_{1})|\,
|\xi_{2}|^s
|\widehat\theta(\xi_{2})|
\\[-2pt]
&\qquad
\times
|\xi_{3}^{(1)}|
|\widehat\theta(\xi_{3})|\,
\prod_{k=4}^{n+2}
|\xi_{k}|
|\widehat\theta(\xi_{k})|
\dd\xi_1
\cdots
\dd\xi_{n+1}.
\end{aligned}
\]

\textit{Case 2:} $i_1=2$, $i_2=3$. Using \eqref{eq:two-largest} and $|\xi_1|^{s-2}\le|\xi_3|^{s-2}$, we have
\[
\begin{aligned}
I_{4,1}^{(n)}
&
\le
(n+1)
\int_{\xi_1+\cdots+\xi_{n+2}=0}
|\xi_{1}^{(1)}|
|\widehat\theta(\xi_{1})|\,
|\xi_{2}^{(1)}|
|\xi_{2}|^{s-1}
|\widehat\theta(\xi_{2})|
\\[-2pt]
&\qquad
\times
|\xi_{3}|^s
|\widehat\theta(\xi_{3})|\,
\prod_{k=4}^{n+2}
|\xi_{k}|
|\widehat\theta(\xi_{k})|
\dd\xi_1
\cdots
\dd\xi_{n+1}.
\end{aligned}
\]
\[
\begin{aligned}
I_{4,2}^{(n)}
&
\le
\int_{\xi_1+\cdots+\xi_{n+2}=0}
|\xi_{1}^{(1)}|
|\widehat\theta(\xi_{1})|\,
|\xi_{2}|^s
|\widehat\theta(\xi_{2})|
\\[-2pt]
&\qquad
\times
|\xi_{3}^{(1)}|
|\xi_{3}|^{s-1}
|\widehat\theta(\xi_{3})|\,
\prod_{k=4}^{n+2}
|\xi_{k}|
|\widehat\theta(\xi_{k})|
\dd\xi_1
\cdots
\dd\xi_{n+1}.
\end{aligned}
\]

\begin{samepage}
\textit{Case 3:} $i_1=2$, $i_2\notin\{1,3\}$. Using \eqref{eq:two-largest} and
\[
|\xi_1|^{s-2}
\le
|\xi_{i_2}|^{s-2},
\qquad
|\xi_3|
\le
|\xi_{i_2}|,
\]
we have
\[
\begin{aligned}
I_{4,1}^{(n)}
&
\le
(n+1)
\int_{\xi_1+\cdots+\xi_{n+2}=0}
|\xi_1^{(1)}|\,
|\widehat\theta(\xi_1)|\,
|\xi_2|^{s-1}
|\xi_2^{(1)}|\,
|\widehat\theta(\xi_2)|
\\
&\qquad
\times
|\xi_3|\,
|\widehat\theta(\xi_3)|\,
|\xi_{i_2}|^s
|\widehat\theta(\xi_{i_2})|
\prod_{\substack{4\le k\le n+2\\k\ne i_2}}
|\xi_k|\,
|\widehat\theta(\xi_k)|
\dd\xi_1
\cdots
\dd\xi_{n+1}
\\
&
\le
(n+1)
\sum_{j=4}^{n+2}
\int_{\xi_1+\cdots+\xi_{n+2}=0}
|\xi_1^{(1)}|\,
|\widehat\theta(\xi_1)|\,
|\xi_2|^{s-1}
|\xi_2^{(1)}|\,
|\widehat\theta(\xi_2)|
\\
&\qquad
\times
|\xi_3|\,
|\widehat\theta(\xi_3)|\,
|\xi_j|^s
|\widehat\theta(\xi_j)|
\prod_{\substack{4\le k\le n+2\\k\ne j}}
|\xi_k|\,
|\widehat\theta(\xi_k)|
\dd\xi_1
\cdots
\dd\xi_{n+1},
\end{aligned}
\]
\[
\begin{aligned}
I_{4,2}^{(n)}
&
\le
\int_{\xi_1+\cdots+\xi_{n+2}=0}
|\xi_1^{(1)}|\,
|\widehat\theta(\xi_1)|\,
|\xi_2|^s
|\widehat\theta(\xi_2)|\,
|\xi_3^{(1)}|\,
|\widehat\theta(\xi_3)|
\\
&\qquad
\times
|\xi_{i_2}|^s
|\widehat\theta(\xi_{i_2})|
\prod_{\substack{4\le k\le n+2\\k\ne i_2}}
|\xi_k|\,
|\widehat\theta(\xi_k)|
\dd\xi_1
\cdots
\dd\xi_{n+1}
\\
&
\le
\sum_{j=4}^{n+2}
\int_{\xi_1+\cdots+\xi_{n+2}=0}
|\xi_1^{(1)}|\,
|\widehat\theta(\xi_1)|\,
|\xi_2|^s
|\widehat\theta(\xi_2)|\,
|\xi_3^{(1)}|\,
|\widehat\theta(\xi_3)|
\\
&\qquad
\times
|\xi_j|^s
|\widehat\theta(\xi_j)|
\prod_{\substack{4\le k\le n+2\\k\ne j}}
|\xi_k|\,
|\widehat\theta(\xi_k)|
\dd\xi_1
\cdots
\dd\xi_{n+1}.
\end{aligned}
\]

\end{samepage}

\textit{Case 4:} $i_1=1$. Since $|\xi_2|,|\xi_3|\le|\xi_1|$, we have
\[
\begin{aligned}
I_{4,1}^{(n)}
&
\le
\int_{\xi_1+\cdots+\xi_{n+2}=0}
|\xi_{1}^{(1)}|
|\xi_{1}|^{s-1}
|\widehat\theta(\xi_{1})|\,
|\xi_{2}^{(1)}|
|\xi_{2}|^{s-1}
|\widehat\theta(\xi_{2})|
\\[-2pt]
&\qquad
\times
\prod_{k=3}^{n+2}
|\xi_{k}|
|\widehat\theta(\xi_{k})|
\dd\xi_1
\cdots
\dd\xi_{n+1}.
\end{aligned}
\]
\[
\begin{aligned}
I_{4,2}^{(n)}
&
\le
\int_{\xi_1+\cdots+\xi_{n+2}=0}
|\xi_{1}^{(1)}|
|\xi_{1}|^{s-1}
|\widehat\theta(\xi_{1})|\,
|\xi_{2}|^s
|\widehat\theta(\xi_{2})|
\\[-2pt]
&\qquad
\times
|\xi_{3}^{(1)}|
|\widehat\theta(\xi_{3})|\,
\prod_{k=4}^{n+2}
|\xi_{k}|
|\widehat\theta(\xi_{k})|
\dd\xi_1
\cdots
\dd\xi_{n+1}.
\end{aligned}
\]

\textit{Case 5:} $i_1=3$. Since $|\xi_1|,|\xi_2|\le|\xi_3|$, we have
\[
\begin{aligned}
I_{4,1}^{(n)}
&
\le
\int_{\xi_1+\cdots+\xi_{n+2}=0}
|\xi_{1}^{(1)}|
|\widehat\theta(\xi_{1})|\,
|\xi_{2}^{(1)}|
|\xi_{2}|^{s-1}
|\widehat\theta(\xi_{2})|
\\[-2pt]
&\qquad
\times
|\xi_{3}|^s
|\widehat\theta(\xi_{3})|\,
\prod_{k=4}^{n+2}
|\xi_{k}|
|\widehat\theta(\xi_{k})|
\dd\xi_1
\cdots
\dd\xi_{n+1}.
\end{aligned}
\]
\[
\begin{aligned}
I_{4,2}^{(n)}
&
\le
\int_{\xi_1+\cdots+\xi_{n+2}=0}
|\xi_{1}^{(1)}|
|\widehat\theta(\xi_{1})|\,
|\xi_{2}|^s
|\widehat\theta(\xi_{2})|
\\[-2pt]
&\qquad
\times
|\xi_{3}^{(1)}|
|\xi_{3}|^{s-1}
|\widehat\theta(\xi_{3})|\,
\prod_{k=4}^{n+2}
|\xi_{k}|
|\widehat\theta(\xi_{k})|
\dd\xi_1
\cdots
\dd\xi_{n+1}.
\end{aligned}
\]

\begin{samepage}
\textit{Case 6:} $i_1\notin\{1,2,3\}$. Using
\[
|\xi_1|^{s-2}
\le
|\xi_{i_1}|^{s-2},
\qquad
|\xi_2|
\le
|\xi_{i_1}|,
\qquad
|\xi_3|
\le
|\xi_{i_1}|,
\]
we have
\[
\begin{aligned}
I_{4,1}^{(n)}
&
\le
\int_{\xi_1+\cdots+\xi_{n+2}=0}
|\xi_1^{(1)}|\,
|\widehat\theta(\xi_1)|\,
|\xi_2^{(1)}|\,
|\xi_2|^{s-1}
|\widehat\theta(\xi_2)|
\\
&\qquad
\times
|\xi_3|\,
|\widehat\theta(\xi_3)|\,
|\xi_{i_1}|^s
|\widehat\theta(\xi_{i_1})|
\prod_{\substack{4\le k\le n+2\\k\ne i_1}}
|\xi_k|\,
|\widehat\theta(\xi_k)|
\dd\xi_1
\cdots
\dd\xi_{n+1}
\\
&
\le
\sum_{j=4}^{n+2}
\int_{\xi_1+\cdots+\xi_{n+2}=0}
|\xi_1^{(1)}|\,
|\widehat\theta(\xi_1)|\,
|\xi_2^{(1)}|\,
|\xi_2|^{s-1}
|\widehat\theta(\xi_2)|
\\
&\qquad
\times
|\xi_3|\,
|\widehat\theta(\xi_3)|\,
|\xi_j|^s
|\widehat\theta(\xi_j)|
\prod_{\substack{4\le k\le n+2\\k\ne j}}
|\xi_k|\,
|\widehat\theta(\xi_k)|
\dd\xi_1
\cdots
\dd\xi_{n+1},
\end{aligned}
\]
\[
\begin{aligned}
I_{4,2}^{(n)}
&
\le
\int_{\xi_1+\cdots+\xi_{n+2}=0}
|\xi_1^{(1)}|\,
|\widehat\theta(\xi_1)|\,
|\xi_2|^s
|\widehat\theta(\xi_2)|\,
|\xi_3^{(1)}|\,
|\widehat\theta(\xi_3)|
\\
&\qquad
\times
|\xi_{i_1}|^s
|\widehat\theta(\xi_{i_1})|
\prod_{\substack{4\le k\le n+2\\k\ne i_1}}
|\xi_k|\,
|\widehat\theta(\xi_k)|
\dd\xi_1
\cdots
\dd\xi_{n+1}
\\
&
\le
\sum_{j=4}^{n+2}
\int_{\xi_1+\cdots+\xi_{n+2}=0}
|\xi_1^{(1)}|\,
|\widehat\theta(\xi_1)|\,
|\xi_2|^s
|\widehat\theta(\xi_2)|\,
|\xi_3^{(1)}|\,
|\widehat\theta(\xi_3)|
\\
&\qquad
\times
|\xi_j|^s
|\widehat\theta(\xi_j)|
\prod_{\substack{4\le k\le n+2\\k\ne j}}
|\xi_k|\,
|\widehat\theta(\xi_k)|
\dd\xi_1
\cdots
\dd\xi_{n+1}.
\end{aligned}
\]

\end{samepage}

Collecting the estimates in Cases 1--6 and then using \eqref{eq:young} and \eqref{eq:L1theta}--\eqref{eq:L2horizontal} yields
\begin{align}
I_{4,1}^{(n)}
&
\le
\bigl[
2(n+1)
+
(n+1)(n-1)
+
2
+
(n-1)
\bigr]
\left(
\sqrt{\frac{\pi}{s-2}}
\norm{\theta}_{H^s}
\right)^n
\norm{\mathcal R_1\theta}_{H^s}^2
\notag
\\
&
=
(n+1)(n+2)
\left(
\sqrt{\frac{\pi}{s-2}}
\norm{\theta}_{H^s}
\right)^n
\norm{\mathcal R_1\theta}_{H^s}^2,
\label{eq:J1bound}
\\
I_{4,2}^{(n)}
&
\le
[4+2(n-1)]
\left(
\sqrt{\frac{\pi}{s-2}}
\norm{\theta}_{H^s}
\right)^n
\norm{\mathcal R_1\theta}_{H^s}^2
\notag
\\
&
=
2(n+1)
\left(
\sqrt{\frac{\pi}{s-2}}
\norm{\theta}_{H^s}
\right)^n
\norm{\mathcal R_1\theta}_{H^s}^2.
\label{eq:J2bound}
\end{align}
Therefore,
\begin{equation}
\label{eq:Jbound}
\begin{aligned}
&
\int_{\xi_1+\cdots+\xi_{n+2}=0}
|\xi_1^{(1)}|^2
|\xi_1|^{s-2}
|\widehat\theta(\xi_1)|
|\xi_2|^s
|\widehat\theta(\xi_2)|
\prod_{k=3}^{n+2}
|\xi_k|
|\widehat\theta(\xi_k)|
\dd\xi_1
\cdots
\dd\xi_{n+1}
\\
&\quad
\le
(n+1)(3n+2)
\left(
\sqrt{\frac{\pi}{s-2}}
\norm{\theta}_{H^s}
\right)^n
\norm{\mathcal R_1\theta}_{H^s}^2
\\
&\quad
\le
3(n+1)^2
\left(
\sqrt{\frac{\pi}{s-2}}
\norm{\theta}_{H^s}
\right)^n
\norm{\mathcal R_1\theta}_{H^s}^2.
\end{aligned}
\end{equation}
Inserting \eqref{eq:Jbound} into \eqref{eq:I4-exp}, we obtain
\begin{equation}
\label{eq:I4}
\begin{aligned}
I_4(T_1,T_2)
&
\le
3
\sum_{n=1}^{\infty}
|c_n|(n+1)^2
\times
\left(
\frac{\norm{\theta}_{L^\infty(T_1,T_2;H^s)}}{2\sqrt{\pi(s-2)}}
\right)^n
\norm{\mathcal R_1\theta}_{L^2(T_1,T_2;H^s)}^2
\\
&
\le
3
\sum_{n=1}^{\infty}
|c_n|(n+1)^2
\alpha^n
\norm{\mathcal R_1\theta}_{L^2(T_1,T_2;H^s)}^2.
\end{aligned}
\end{equation}

\paragraph{Estimates for \texorpdfstring{$I_5$}{I5}.}
Expanding $\theta_s$ and using \eqref{eq:second-diff}, we obtain
\begin{align}
I_5(T_1,T_2)
&\leq
\sum_{n=0}^{\infty}|c_n|
\int_{T_1}^{T_2}
\left|
\int_{\mathbb R^2}
R_v\,\Lambda_2^s\theta\,(\partial_2\theta)^n\,dx
\right|\,d\tau
\notag\\
&\leq
s(s-1)
2^{s-2}
\sum_{n=0}^{\infty}
|c_n|(2\pi)^{-(n+1)}
\int_{T_1}^{T_2}
\int_{\xi_1+\cdots+\xi_{n+3}=0}
\notag
\\*[-2pt]
&\quad
\bigl(
|\xi_1^{(2)}|^2
|\xi_2^{(2)}|^{s-2}
+
|\xi_1^{(2)}|^s
\bigr)
|\widehat{u^2}(\xi_1)|\,
|\xi_2|
|\widehat\theta(\xi_2)|
\notag
\\*[-2pt]
&\qquad
\times
|\xi_3|^s
|\widehat\theta(\xi_3)|
\prod_{k=4}^{n+3}
|\xi_k|
|\widehat\theta(\xi_k)|
\dd\xi_1
\cdots
\dd\xi_{n+2}
\dd\tau
\notag
\\
&
\le
s(s-1)
2^{s-2}
\sum_{n=0}^{\infty}
|c_n|(2\pi)^{-(n+1)}
\int_{T_1}^{T_2}
\int_{\xi_1+\cdots+\xi_{n+3}=0}
|\xi_1^{(1)}|^2
|\widehat\theta(\xi_1)|
\notag
\\*[-2pt]
&\qquad
\times
|\xi_2|^{s-1}
|\widehat\theta(\xi_2)|
|\xi_3|^s
|\widehat\theta(\xi_3)|
\prod_{k=4}^{n+3}
|\xi_k|
|\widehat\theta(\xi_k)|
\dd\xi_1
\cdots
\dd\xi_{n+2}
\dd\tau
\notag
\\
&\quad
+
s(s-1)
2^{s-2}
\sum_{n=0}^{\infty}
|c_n|(2\pi)^{-(n+1)}
\int_{T_1}^{T_2}
\int_{\xi_1+\cdots+\xi_{n+3}=0}
|\xi_1^{(1)}|^2
|\xi_1|^{s-2}
|\widehat\theta(\xi_1)|
\notag
\\*[-2pt]
&\qquad
\times
|\xi_2|
|\widehat\theta(\xi_2)|
|\xi_3|^s
|\widehat\theta(\xi_3)|
\prod_{k=4}^{n+3}
|\xi_k|
|\widehat\theta(\xi_k)|
\dd\xi_1
\cdots
\dd\xi_{n+2}
\dd\tau.
\label{eq:I5-exp}
\end{align}
From \eqref{eq:Jbound}, we have
\begin{equation}
\label{eq:Vbound}
\begin{aligned}
&
\int_{\xi_1+\cdots+\xi_{n+3}=0}
|\xi_1^{(1)}|^2
|\xi_1|^{s-2}
|\widehat\theta(\xi_1)|
|\xi_2|
|\widehat\theta(\xi_2)|\,
|\xi_3|^s
|\widehat\theta(\xi_3)|
\\[-2pt]
&\qquad
\times
\prod_{k=4}^{n+3}
|\xi_k|
|\widehat\theta(\xi_k)|
\dd\xi_1
\cdots
\dd\xi_{n+2}
\\
&\quad
\le
3(n+2)^2
\left(
\sqrt{\frac{\pi}{s-2}}
\norm{\theta}_{H^s}
\right)^{n+1}
\norm{\mathcal R_1\theta}_{H^s}^2.
\end{aligned}
\end{equation}
For the second integral, using
$|\xi_1^{(1)}|\le\sum_{j=2}^{n+3}|\xi_j^{(1)}|$ and symmetry, we obtain
\begin{align}
&
\int_{\xi_1+\cdots+\xi_{n+3}=0}
|\xi_1^{(1)}|^2
|\widehat\theta(\xi_1)|
|\xi_2|^{s-1}
|\widehat\theta(\xi_2)|
|\xi_3|^s
|\widehat\theta(\xi_3)|
\prod_{k=4}^{n+3}
|\xi_k|
|\widehat\theta(\xi_k)|
\dd\xi_1
\cdots
\dd\xi_{n+2}
\notag
\\
&\quad
\le
\int_{\xi_1+\cdots+\xi_{n+3}=0}
|\xi_1^{(1)}|
|\widehat\theta(\xi_1)|
|\xi_2^{(1)}|
|\xi_2|^{s-1}
|\widehat\theta(\xi_2)|
\notag
\\*[-2pt]
&\qquad\qquad
\times
|\xi_3|^s
|\widehat\theta(\xi_3)|
\prod_{k=4}^{n+3}
|\xi_k|
|\widehat\theta(\xi_k)|
\dd\xi_1
\cdots
\dd\xi_{n+2}
\notag
\\
&\qquad
+
\int_{\xi_1+\cdots+\xi_{n+3}=0}
|\xi_1^{(1)}|
|\widehat\theta(\xi_1)|
|\xi_2|^{s-1}
|\widehat\theta(\xi_2)|
\notag
\\*[-2pt]
&\qquad\qquad
\times
|\xi_3^{(1)}|
|\xi_3|^s
|\widehat\theta(\xi_3)|
\prod_{k=4}^{n+3}
|\xi_k|
|\widehat\theta(\xi_k)|
\dd\xi_1
\cdots
\dd\xi_{n+2}
\notag
\\
&\qquad
+
n
\int_{\xi_1+\cdots+\xi_{n+3}=0}
|\xi_1^{(1)}|
|\widehat\theta(\xi_1)|
|\xi_2|^{s-1}
|\widehat\theta(\xi_2)|
|\xi_3|^s
|\widehat\theta(\xi_3)|
\notag
\\*[-2pt]
&\qquad\qquad
\times
|\xi_4^{(1)}|
|\xi_4|
|\widehat\theta(\xi_4)|
\prod_{k=5}^{n+3}
|\xi_k|
|\widehat\theta(\xi_k)|
\dd\xi_1
\cdots
\dd\xi_{n+2}
\notag
\\
&\quad
=:
I_{5,1}^{(n)}
+
I_{5,2}^{(n)}
+
nI_{5,3}^{(n)}.
\label{eq:W123}
\end{align}
Let
$|\xi_{i_1}|$ and $|\xi_{i_2}|$ be the largest and second-largest of
$|\xi_1|,\ldots,|\xi_{n+3}|$. Then
$|\xi_{i_1}|\le(n+2)|\xi_{i_2}|$. We consider the following cases.

\textit{Case 1:} $i_1=3$ and $i_2\in\{1,2,4\}$.
If $i_2=1$, we have
\begin{equation*}
\begin{aligned}
I_{5,1}^{(n)}
&
\le
\int_{\xi_1+\cdots+\xi_{n+3}=0}
|\xi_{1}^{(1)}|
|\xi_{1}|^{s-1}
|\widehat\theta(\xi_{1})|\,
|\xi_{2}^{(1)}|
|\widehat\theta(\xi_{2})|
\\[-2pt]
&\qquad
\times
|\xi_{3}|^s
|\widehat\theta(\xi_{3})|\,
|\xi_{4}|
|\widehat\theta(\xi_{4})|\,
\prod_{k=5}^{n+3}
|\xi_{k}|
|\widehat\theta(\xi_{k})|
\dd\xi_1
\cdots
\dd\xi_{n+2}.
\end{aligned}
\end{equation*}
\begin{equation*}
\begin{aligned}
I_{5,2}^{(n)}
&
\le
(n+2)
\int_{\xi_1+\cdots+\xi_{n+3}=0}
|\xi_{1}^{(1)}|
|\xi_{1}|^{s-1}
|\widehat\theta(\xi_{1})|\,
|\xi_{2}|
|\widehat\theta(\xi_{2})|
\\[-2pt]
&\qquad
\times
|\xi_{3}^{(1)}|
|\xi_{3}|^{s-1}
|\widehat\theta(\xi_{3})|\,
|\xi_{4}|
|\widehat\theta(\xi_{4})|\,
\prod_{k=5}^{n+3}
|\xi_{k}|
|\widehat\theta(\xi_{k})|
\dd\xi_1
\cdots
\dd\xi_{n+2}.
\end{aligned}
\end{equation*}
\begin{equation*}
\begin{aligned}
I_{5,3}^{(n)}
&
\le
\int_{\xi_1+\cdots+\xi_{n+3}=0}
|\xi_{1}^{(1)}|
|\xi_{1}|^{s-1}
|\widehat\theta(\xi_{1})|\,
|\xi_{2}|
|\widehat\theta(\xi_{2})|
\\[-2pt]
&\qquad
\times
|\xi_{3}|^s
|\widehat\theta(\xi_{3})|\,
|\xi_{4}^{(1)}|
|\widehat\theta(\xi_{4})|\,
\prod_{k=5}^{n+3}
|\xi_{k}|
|\widehat\theta(\xi_{k})|
\dd\xi_1
\cdots
\dd\xi_{n+2}.
\end{aligned}
\end{equation*}
If $i_2=2$, we have
\begin{equation*}
\begin{aligned}
I_{5,1}^{(n)}
&
\le
\int_{\xi_1+\cdots+\xi_{n+3}=0}
|\xi_{1}^{(1)}|
|\widehat\theta(\xi_{1})|\,
|\xi_{2}^{(1)}|
|\xi_{2}|^{s-1}
|\widehat\theta(\xi_{2})|
\\[-2pt]
&\qquad
\times
|\xi_{3}|^s
|\widehat\theta(\xi_{3})|\,
|\xi_{4}|
|\widehat\theta(\xi_{4})|\,
\prod_{k=5}^{n+3}
|\xi_{k}|
|\widehat\theta(\xi_{k})|
\dd\xi_1
\cdots
\dd\xi_{n+2}.
\end{aligned}
\end{equation*}
\begin{equation*}
\begin{aligned}
I_{5,2}^{(n)}
&
\le
(n+2)
\int_{\xi_1+\cdots+\xi_{n+3}=0}
|\xi_{1}^{(1)}|
|\widehat\theta(\xi_{1})|\,
|\xi_{2}|^s
|\widehat\theta(\xi_{2})|
\\[-2pt]
&\qquad
\times
|\xi_{3}^{(1)}|
|\xi_{3}|^{s-1}
|\widehat\theta(\xi_{3})|\,
|\xi_{4}|
|\widehat\theta(\xi_{4})|\,
\prod_{k=5}^{n+3}
|\xi_{k}|
|\widehat\theta(\xi_{k})|
\dd\xi_1
\cdots
\dd\xi_{n+2}.
\end{aligned}
\end{equation*}
\begin{equation*}
\begin{aligned}
I_{5,3}^{(n)}
&
\le
\int_{\xi_1+\cdots+\xi_{n+3}=0}
|\xi_{1}^{(1)}|
|\widehat\theta(\xi_{1})|\,
|\xi_{2}|^s
|\widehat\theta(\xi_{2})|
\\[-2pt]
&\qquad
\times
|\xi_{3}|^s
|\widehat\theta(\xi_{3})|\,
|\xi_{4}^{(1)}|
|\widehat\theta(\xi_{4})|\,
\prod_{k=5}^{n+3}
|\xi_{k}|
|\widehat\theta(\xi_{k})|
\dd\xi_1
\cdots
\dd\xi_{n+2}.
\end{aligned}
\end{equation*}
If $i_2=4$, we have
\begin{equation*}
\begin{aligned}
I_{5,1}^{(n)}
&
\le
\int_{\xi_1+\cdots+\xi_{n+3}=0}
|\xi_{1}^{(1)}|
|\widehat\theta(\xi_{1})|\,
|\xi_{2}^{(1)}|
|\widehat\theta(\xi_{2})|
\\[-2pt]
&\qquad
\times
|\xi_{3}|^s
|\widehat\theta(\xi_{3})|\,
|\xi_{4}|^s
|\widehat\theta(\xi_{4})|\,
\prod_{k=5}^{n+3}
|\xi_{k}|
|\widehat\theta(\xi_{k})|
\dd\xi_1
\cdots
\dd\xi_{n+2}.
\end{aligned}
\end{equation*}
\begin{equation*}
\begin{aligned}
I_{5,2}^{(n)}
&
\le
(n+2)
\int_{\xi_1+\cdots+\xi_{n+3}=0}
|\xi_{1}^{(1)}|
|\widehat\theta(\xi_{1})|\,
|\xi_{2}|
|\widehat\theta(\xi_{2})|
\\[-2pt]
&\qquad
\times
|\xi_{3}^{(1)}|
|\xi_{3}|^{s-1}
|\widehat\theta(\xi_{3})|\,
|\xi_{4}|^s
|\widehat\theta(\xi_{4})|\,
\prod_{k=5}^{n+3}
|\xi_{k}|
|\widehat\theta(\xi_{k})|
\dd\xi_1
\cdots
\dd\xi_{n+2}.
\end{aligned}
\end{equation*}
\begin{equation*}
\begin{aligned}
I_{5,3}^{(n)}
&
\le
\int_{\xi_1+\cdots+\xi_{n+3}=0}
|\xi_{1}^{(1)}|
|\widehat\theta(\xi_{1})|\,
|\xi_{2}|
|\widehat\theta(\xi_{2})|
\\[-2pt]
&\qquad
\times
|\xi_{3}|^s
|\widehat\theta(\xi_{3})|\,
|\xi_{4}^{(1)}|
|\xi_{4}|^{s-1}
|\widehat\theta(\xi_{4})|\,
\prod_{k=5}^{n+3}
|\xi_{k}|
|\widehat\theta(\xi_{k})|
\dd\xi_1
\cdots
\dd\xi_{n+2}.
\end{aligned}
\end{equation*}

\textit{Case 2:} $i_1=3$ and $5\le i_2\le n+3$.
We have
\begin{equation*}
\begin{aligned}
I_{5,1}^{(n)}
&
\le
\int_{\xi_1+\cdots+\xi_{n+3}=0}
|\xi_{1}^{(1)}|
|\widehat\theta(\xi_{1})|\,
|\xi_{2}^{(1)}|
|\widehat\theta(\xi_{2})|
\\[-2pt]
&\qquad
\times
|\xi_{3}|^s
|\widehat\theta(\xi_{3})|\,
|\xi_{4}|
|\widehat\theta(\xi_{4})|\,
|\xi_{i_2}|^s
|\widehat\theta(\xi_{i_2})|\,
\prod_{\substack{k\in\{5,\ldots,n+3\}\\k\ne i_2}}
|\xi_{k}|
|\widehat\theta(\xi_{k})|
\dd\xi_1
\cdots
\dd\xi_{n+2}
\\
&
\le
\sum_{j=5}^{n+3}
\int_{\xi_1+\cdots+\xi_{n+3}=0}
|\xi_{1}^{(1)}|
|\widehat\theta(\xi_{1})|\,
|\xi_{2}^{(1)}|
|\widehat\theta(\xi_{2})|
\\[-2pt]
&\qquad
\times
|\xi_{3}|^s
|\widehat\theta(\xi_{3})|\,
|\xi_{4}|
|\widehat\theta(\xi_{4})|\,
|\xi_j|^s
|\widehat\theta(\xi_j)|\,
\prod_{\substack{k\in\{5,\ldots,n+3\}\\k\ne j}}
|\xi_{k}|
|\widehat\theta(\xi_{k})|
\dd\xi_1
\cdots
\dd\xi_{n+2}.
\end{aligned}
\end{equation*}
\begin{equation*}
\begin{aligned}
I_{5,2}^{(n)}
&
\le
(n+2)
\int_{\xi_1+\cdots+\xi_{n+3}=0}
|\xi_{1}^{(1)}|
|\widehat\theta(\xi_{1})|\,
|\xi_{2}|
|\widehat\theta(\xi_{2})|
\\[-2pt]
&\qquad
\times
|\xi_{3}^{(1)}|
|\xi_{3}|^{s-1}
|\widehat\theta(\xi_{3})|\,
|\xi_{4}|
|\widehat\theta(\xi_{4})|\,
|\xi_{i_2}|^s
|\widehat\theta(\xi_{i_2})|\,
\prod_{\substack{k\in\{5,\ldots,n+3\}\\k\ne i_2}}
|\xi_{k}|
|\widehat\theta(\xi_{k})|
\dd\xi_1
\cdots
\dd\xi_{n+2}
\\
&
\le
(n+2)
\sum_{j=5}^{n+3}
\int_{\xi_1+\cdots+\xi_{n+3}=0}
|\xi_{1}^{(1)}|
|\widehat\theta(\xi_{1})|\,
|\xi_{2}|
|\widehat\theta(\xi_{2})|
\\[-2pt]
&\qquad
\times
|\xi_{3}^{(1)}|
|\xi_{3}|^{s-1}
|\widehat\theta(\xi_{3})|\,
|\xi_{4}|
|\widehat\theta(\xi_{4})|\,
|\xi_j|^s
|\widehat\theta(\xi_j)|\,
\prod_{\substack{k\in\{5,\ldots,n+3\}\\k\ne j}}
|\xi_{k}|
|\widehat\theta(\xi_{k})|
\dd\xi_1
\cdots
\dd\xi_{n+2}.
\end{aligned}
\end{equation*}
\begin{equation*}
\begin{aligned}
I_{5,3}^{(n)}
&
\le
\int_{\xi_1+\cdots+\xi_{n+3}=0}
|\xi_{1}^{(1)}|
|\widehat\theta(\xi_{1})|\,
|\xi_{2}|
|\widehat\theta(\xi_{2})|
\\[-2pt]
&\qquad
\times
|\xi_{3}|^s
|\widehat\theta(\xi_{3})|\,
|\xi_{4}^{(1)}|
|\widehat\theta(\xi_{4})|\,
|\xi_{i_2}|^s
|\widehat\theta(\xi_{i_2})|\,
\prod_{\substack{k\in\{5,\ldots,n+3\}\\k\ne i_2}}
|\xi_{k}|
|\widehat\theta(\xi_{k})|
\dd\xi_1
\cdots
\dd\xi_{n+2}
\\
&
\le
\sum_{j=5}^{n+3}
\int_{\xi_1+\cdots+\xi_{n+3}=0}
|\xi_{1}^{(1)}|
|\widehat\theta(\xi_{1})|\,
|\xi_{2}|
|\widehat\theta(\xi_{2})|
\\[-2pt]
&\qquad
\times
|\xi_{3}|^s
|\widehat\theta(\xi_{3})|\,
|\xi_{4}^{(1)}|
|\widehat\theta(\xi_{4})|\,
|\xi_j|^s
|\widehat\theta(\xi_j)|\,
\prod_{\substack{k\in\{5,\ldots,n+3\}\\k\ne j}}
|\xi_{k}|
|\widehat\theta(\xi_{k})|
\dd\xi_1
\cdots
\dd\xi_{n+2}.
\end{aligned}
\end{equation*}

\textit{Case 3:} $i_1\in\{1,2,4\}$.
If $i_1=1$, we have
\begin{equation*}
\begin{aligned}
I_{5,1}^{(n)}
&
\le
\int_{\xi_1+\cdots+\xi_{n+3}=0}
|\xi_{1}^{(1)}|
|\xi_{1}|^{s-1}
|\widehat\theta(\xi_{1})|\,
|\xi_{2}^{(1)}|
|\widehat\theta(\xi_{2})|
\\[-2pt]
&\qquad
\times
|\xi_{3}|^s
|\widehat\theta(\xi_{3})|\,
|\xi_{4}|
|\widehat\theta(\xi_{4})|\,
\prod_{k=5}^{n+3}
|\xi_{k}|
|\widehat\theta(\xi_{k})|
\dd\xi_1
\cdots
\dd\xi_{n+2}.
\end{aligned}
\end{equation*}
\begin{equation*}
\begin{aligned}
I_{5,2}^{(n)}
&
\le
\int_{\xi_1+\cdots+\xi_{n+3}=0}
|\xi_{1}^{(1)}|
|\xi_{1}|^{s-1}
|\widehat\theta(\xi_{1})|\,
|\xi_{2}|
|\widehat\theta(\xi_{2})|
\\[-2pt]
&\qquad
\times
|\xi_{3}^{(1)}|
|\xi_{3}|^{s-1}
|\widehat\theta(\xi_{3})|\,
|\xi_{4}|
|\widehat\theta(\xi_{4})|\,
\prod_{k=5}^{n+3}
|\xi_{k}|
|\widehat\theta(\xi_{k})|
\dd\xi_1
\cdots
\dd\xi_{n+2}.
\end{aligned}
\end{equation*}
\begin{equation*}
\begin{aligned}
I_{5,3}^{(n)}
&
\le
\int_{\xi_1+\cdots+\xi_{n+3}=0}
|\xi_{1}^{(1)}|
|\xi_{1}|^{s-1}
|\widehat\theta(\xi_{1})|\,
|\xi_{2}|
|\widehat\theta(\xi_{2})|
\\[-2pt]
&\qquad
\times
|\xi_{3}|^s
|\widehat\theta(\xi_{3})|\,
|\xi_{4}^{(1)}|
|\widehat\theta(\xi_{4})|\,
\prod_{k=5}^{n+3}
|\xi_{k}|
|\widehat\theta(\xi_{k})|
\dd\xi_1
\cdots
\dd\xi_{n+2}.
\end{aligned}
\end{equation*}
If $i_1=2$, we have
\begin{equation*}
\begin{aligned}
I_{5,1}^{(n)}
&
\le
\int_{\xi_1+\cdots+\xi_{n+3}=0}
|\xi_{1}^{(1)}|
|\widehat\theta(\xi_{1})|\,
|\xi_{2}^{(1)}|
|\xi_{2}|^{s-1}
|\widehat\theta(\xi_{2})|
\\[-2pt]
&\qquad
\times
|\xi_{3}|^s
|\widehat\theta(\xi_{3})|\,
|\xi_{4}|
|\widehat\theta(\xi_{4})|\,
\prod_{k=5}^{n+3}
|\xi_{k}|
|\widehat\theta(\xi_{k})|
\dd\xi_1
\cdots
\dd\xi_{n+2}.
\end{aligned}
\end{equation*}
\begin{equation*}
\begin{aligned}
I_{5,2}^{(n)}
&
\le
\int_{\xi_1+\cdots+\xi_{n+3}=0}
|\xi_{1}^{(1)}|
|\widehat\theta(\xi_{1})|\,
|\xi_{2}|^s
|\widehat\theta(\xi_{2})|
\\[-2pt]
&\qquad
\times
|\xi_{3}^{(1)}|
|\xi_{3}|^{s-1}
|\widehat\theta(\xi_{3})|\,
|\xi_{4}|
|\widehat\theta(\xi_{4})|\,
\prod_{k=5}^{n+3}
|\xi_{k}|
|\widehat\theta(\xi_{k})|
\dd\xi_1
\cdots
\dd\xi_{n+2}.
\end{aligned}
\end{equation*}
\begin{equation*}
\begin{aligned}
I_{5,3}^{(n)}
&
\le
\int_{\xi_1+\cdots+\xi_{n+3}=0}
|\xi_{1}^{(1)}|
|\widehat\theta(\xi_{1})|\,
|\xi_{2}|^s
|\widehat\theta(\xi_{2})|
\\[-2pt]
&\qquad
\times
|\xi_{3}|^s
|\widehat\theta(\xi_{3})|\,
|\xi_{4}^{(1)}|
|\widehat\theta(\xi_{4})|\,
\prod_{k=5}^{n+3}
|\xi_{k}|
|\widehat\theta(\xi_{k})|
\dd\xi_1
\cdots
\dd\xi_{n+2}.
\end{aligned}
\end{equation*}
If $i_1=4$, we have
\begin{equation*}
\begin{aligned}
I_{5,1}^{(n)}
&
\le
\int_{\xi_1+\cdots+\xi_{n+3}=0}
|\xi_{1}^{(1)}|
|\widehat\theta(\xi_{1})|\,
|\xi_{2}^{(1)}|
|\widehat\theta(\xi_{2})|
\\[-2pt]
&\qquad
\times
|\xi_{3}|^s
|\widehat\theta(\xi_{3})|\,
|\xi_{4}|^s
|\widehat\theta(\xi_{4})|\,
\prod_{k=5}^{n+3}
|\xi_{k}|
|\widehat\theta(\xi_{k})|
\dd\xi_1
\cdots
\dd\xi_{n+2}.
\end{aligned}
\end{equation*}
\begin{equation*}
\begin{aligned}
I_{5,2}^{(n)}
&
\le
\int_{\xi_1+\cdots+\xi_{n+3}=0}
|\xi_{1}^{(1)}|
|\widehat\theta(\xi_{1})|\,
|\xi_{2}|
|\widehat\theta(\xi_{2})|
\\[-2pt]
&\qquad
\times
|\xi_{3}^{(1)}|
|\xi_{3}|^{s-1}
|\widehat\theta(\xi_{3})|\,
|\xi_{4}|^s
|\widehat\theta(\xi_{4})|\,
\prod_{k=5}^{n+3}
|\xi_{k}|
|\widehat\theta(\xi_{k})|
\dd\xi_1
\cdots
\dd\xi_{n+2}.
\end{aligned}
\end{equation*}
\begin{equation*}
\begin{aligned}
I_{5,3}^{(n)}
&
\le
\int_{\xi_1+\cdots+\xi_{n+3}=0}
|\xi_{1}^{(1)}|
|\widehat\theta(\xi_{1})|\,
|\xi_{2}|
|\widehat\theta(\xi_{2})|
\\[-2pt]
&\qquad
\times
|\xi_{3}|^s
|\widehat\theta(\xi_{3})|\,
|\xi_{4}^{(1)}|
|\xi_{4}|^{s-1}
|\widehat\theta(\xi_{4})|\,
\prod_{k=5}^{n+3}
|\xi_{k}|
|\widehat\theta(\xi_{k})|
\dd\xi_1
\cdots
\dd\xi_{n+2}.
\end{aligned}
\end{equation*}

\textit{Case 4:} $5\le i_1\le n+3$.
We have
\begin{equation*}
\begin{aligned}
I_{5,1}^{(n)}
&
\le
\int_{\xi_1+\cdots+\xi_{n+3}=0}
|\xi_{1}^{(1)}|
|\widehat\theta(\xi_{1})|\,
|\xi_{2}^{(1)}|
|\widehat\theta(\xi_{2})|
\\[-2pt]
&\qquad
\times
|\xi_{3}|^s
|\widehat\theta(\xi_{3})|\,
|\xi_{4}|
|\widehat\theta(\xi_{4})|\,
|\xi_{i_1}|^s
|\widehat\theta(\xi_{i_1})|\,
\prod_{\substack{k\in\{5,\ldots,n+3\}\\k\ne i_1}}
|\xi_{k}|
|\widehat\theta(\xi_{k})|
\dd\xi_1
\cdots
\dd\xi_{n+2}
\\
&
\le
\sum_{j=5}^{n+3}
\int_{\xi_1+\cdots+\xi_{n+3}=0}
|\xi_{1}^{(1)}|
|\widehat\theta(\xi_{1})|\,
|\xi_{2}^{(1)}|
|\widehat\theta(\xi_{2})|
\\[-2pt]
&\qquad
\times
|\xi_{3}|^s
|\widehat\theta(\xi_{3})|\,
|\xi_{4}|
|\widehat\theta(\xi_{4})|\,
|\xi_j|^s
|\widehat\theta(\xi_j)|\,
\prod_{\substack{k\in\{5,\ldots,n+3\}\\k\ne j}}
|\xi_{k}|
|\widehat\theta(\xi_{k})|
\dd\xi_1
\cdots
\dd\xi_{n+2}.
\end{aligned}
\end{equation*}
\begin{equation*}
\begin{aligned}
I_{5,2}^{(n)}
&
\le
\int_{\xi_1+\cdots+\xi_{n+3}=0}
|\xi_{1}^{(1)}|
|\widehat\theta(\xi_{1})|\,
|\xi_{2}|
|\widehat\theta(\xi_{2})|
\\[-2pt]
&\qquad
\times
|\xi_{3}^{(1)}|
|\xi_{3}|^{s-1}
|\widehat\theta(\xi_{3})|\,
|\xi_{4}|
|\widehat\theta(\xi_{4})|\,
|\xi_{i_1}|^s
|\widehat\theta(\xi_{i_1})|\,
\prod_{\substack{k\in\{5,\ldots,n+3\}\\k\ne i_1}}
|\xi_{k}|
|\widehat\theta(\xi_{k})|
\dd\xi_1
\cdots
\dd\xi_{n+2}
\\
&
\le
\sum_{j=5}^{n+3}
\int_{\xi_1+\cdots+\xi_{n+3}=0}
|\xi_{1}^{(1)}|
|\widehat\theta(\xi_{1})|\,
|\xi_{2}|
|\widehat\theta(\xi_{2})|
\\[-2pt]
&\qquad
\times
|\xi_{3}^{(1)}|
|\xi_{3}|^{s-1}
|\widehat\theta(\xi_{3})|\,
|\xi_{4}|
|\widehat\theta(\xi_{4})|\,
|\xi_j|^s
|\widehat\theta(\xi_j)|\,
\prod_{\substack{k\in\{5,\ldots,n+3\}\\k\ne j}}
|\xi_{k}|
|\widehat\theta(\xi_{k})|
\dd\xi_1
\cdots
\dd\xi_{n+2}.
\end{aligned}
\end{equation*}
\begin{equation*}
\begin{aligned}
I_{5,3}^{(n)}
&
\le
\int_{\xi_1+\cdots+\xi_{n+3}=0}
|\xi_{1}^{(1)}|
|\widehat\theta(\xi_{1})|\,
|\xi_{2}|
|\widehat\theta(\xi_{2})|
\\[-2pt]
&\qquad
\times
|\xi_{3}|^s
|\widehat\theta(\xi_{3})|\,
|\xi_{4}^{(1)}|
|\widehat\theta(\xi_{4})|\,
|\xi_{i_1}|^s
|\widehat\theta(\xi_{i_1})|\,
\prod_{\substack{k\in\{5,\ldots,n+3\}\\k\ne i_1}}
|\xi_{k}|
|\widehat\theta(\xi_{k})|
\dd\xi_1
\cdots
\dd\xi_{n+2}
\\
&
\le
\sum_{j=5}^{n+3}
\int_{\xi_1+\cdots+\xi_{n+3}=0}
|\xi_{1}^{(1)}|
|\widehat\theta(\xi_{1})|\,
|\xi_{2}|
|\widehat\theta(\xi_{2})|
\\[-2pt]
&\qquad
\times
|\xi_{3}|^s
|\widehat\theta(\xi_{3})|\,
|\xi_{4}^{(1)}|
|\widehat\theta(\xi_{4})|\,
|\xi_j|^s
|\widehat\theta(\xi_j)|\,
\prod_{\substack{k\in\{5,\ldots,n+3\}\\k\ne j}}
|\xi_{k}|
|\widehat\theta(\xi_{k})|
\dd\xi_1
\cdots
\dd\xi_{n+2}.
\end{aligned}
\end{equation*}

Collecting the estimates in Cases 1--4 and then using \eqref{eq:young} and \eqref{eq:L1theta}--\eqref{eq:L2horizontal}, we obtain
\begin{align}
I_{5,1}^{(n)}
&
\le
[6+2(n-1)]
\left(
\sqrt{\frac{\pi}{s-2}}
\norm{\theta}_{H^s}
\right)^{n+1}
\norm{\mathcal R_1\theta}_{H^s}^2
\notag
\\
&
=
2(n+2)
\left(
\sqrt{\frac{\pi}{s-2}}
\norm{\theta}_{H^s}
\right)^{n+1}
\norm{\mathcal R_1\theta}_{H^s}^2,
\label{eq:W1bound}
\\
I_{5,2}^{(n)}
&
\le
[3(n+2)+(n+2)(n-1)+3+(n-1)]
\notag
\\*[-2pt]
&\qquad
\times
\left(
\sqrt{\frac{\pi}{s-2}}
\norm{\theta}_{H^s}
\right)^{n+1}
\norm{\mathcal R_1\theta}_{H^s}^2
\notag
\\
&
\le
2(n+2)^2
\left(
\sqrt{\frac{\pi}{s-2}}
\norm{\theta}_{H^s}
\right)^{n+1}
\norm{\mathcal R_1\theta}_{H^s}^2,
\label{eq:W2bound}
\\
I_{5,3}^{(n)}
&
\le
[6+2(n-1)]
\left(
\sqrt{\frac{\pi}{s-2}}
\norm{\theta}_{H^s}
\right)^{n+1}
\norm{\mathcal R_1\theta}_{H^s}^2
\notag
\\
&
=
2(n+2)
\left(
\sqrt{\frac{\pi}{s-2}}
\norm{\theta}_{H^s}
\right)^{n+1}
\norm{\mathcal R_1\theta}_{H^s}^2.
\label{eq:W3bound}
\end{align}
Inserting \eqref{eq:W1bound}--\eqref{eq:W3bound} into \eqref{eq:W123}, we obtain
\begin{equation}
\label{eq:Wbound}
\begin{aligned}
&
\int_{\xi_1+\cdots+\xi_{n+3}=0}
|\xi_1^{(1)}|^2
|\widehat\theta(\xi_1)|
|\xi_2|^{s-1}
|\widehat\theta(\xi_2)|
|\xi_3|^s
|\widehat\theta(\xi_3)|
\\[-2pt]
&\qquad
\times
\prod_{k=4}^{n+3}
|\xi_k|
|\widehat\theta(\xi_k)|
\dd\xi_1
\cdots
\dd\xi_{n+2}
\\
&\quad
\le
8(n+2)^2
\left(
\sqrt{\frac{\pi}{s-2}}
\norm{\theta}_{H^s}
\right)^{n+1}
\norm{\mathcal R_1\theta}_{H^s}^2.
\end{aligned}
\end{equation}
Inserting \eqref{eq:Vbound} and \eqref{eq:Wbound} into \eqref{eq:I5-exp} gives
\begin{align}
I_5(T_1,T_2)
&
\le
11s(s-1)
2^{s-2}
\sum_{n=0}^{\infty}
|c_n|(n+2)^2
\notag
\\*
&\qquad
\times
\left(
\frac{\norm{\theta}_{L^\infty(T_1,T_2;H^s)}}{2\sqrt{\pi(s-2)}}
\right)^{n+1}
\norm{\mathcal R_1\theta}_{L^2(T_1,T_2;H^s)}^2
\notag
\\
&
=
11s(s-1)
2^{s-2}
\sum_{n=1}^{\infty}
|c_{n-1}|(n+1)^2
\notag
\\*
&\qquad
\times
\left(
\frac{\norm{\theta}_{L^\infty(T_1,T_2;H^s)}}{2\sqrt{\pi(s-2)}}
\right)^n
\norm{\mathcal R_1\theta}_{L^2(T_1,T_2;H^s)}^2
\notag
\\
&
\le
11s(s-1)
2^{s-2}
\sum_{n=1}^{\infty}
|c_{n-1}|(n+1)^2
\alpha^n
\norm{\mathcal R_1\theta}_{L^2(T_1,T_2;H^s)}^2.
\label{eq:I5}
\end{align}

\paragraph{Estimates for \texorpdfstring{$I_6$}{I6}.}
Expanding $\theta_s$ and using \eqref{eq:first-diff}, we obtain
\begin{align}
I_6(T_1,T_2)
&\leq
\sum_{n=0}^{\infty}|c_n|
\int_{T_1}^{T_2}
\left|
\int_{\mathbb R^2}
R_h\,\Lambda_2^s\theta\,(\partial_2\theta)^n\,dx
\right|\,d\tau
\notag\\
&\leq
s2^{s-1}
\sum_{n=0}^{\infty}
|c_n|(2\pi)^{-(n+1)}
\int_{T_1}^{T_2}
\int_{\xi_1+\cdots+\xi_{n+3}=0}
\notag
\\*[-2pt]
&\quad
\bigl(
|\xi_1^{(2)}|
|\xi_2^{(2)}|^{s-1}
+
|\xi_1^{(2)}|^s
\bigr)
|\widehat{u^1}(\xi_1)|\,
|\xi_2^{(1)}|
|\widehat\theta(\xi_2)|
\notag
\\*[-2pt]
&\qquad
\times
|\xi_3^{(2)}|^s
|\widehat\theta(\xi_3)|
\prod_{k=4}^{n+3}
|\xi_k|
|\widehat\theta(\xi_k)|
\dd\xi_1
\cdots
\dd\xi_{n+2}
\dd\tau
\notag
\\
&
\le
s2^{s-1}
\sum_{n=0}^{\infty}
|c_n|(2\pi)^{-(n+1)}
\int_{T_1}^{T_2}
\int_{\xi_1+\cdots+\xi_{n+3}=0}
|\xi_1^{(1)}|
|\widehat\theta(\xi_1)|
\notag
\\*[-2pt]
&\qquad
\times
|\xi_2^{(1)}|
|\xi_2|^{s-1}
|\widehat\theta(\xi_2)|
|\xi_3|^s
|\widehat\theta(\xi_3)|
\prod_{k=4}^{n+3}
|\xi_k|
|\widehat\theta(\xi_k)|
\dd\xi_1
\cdots
\dd\xi_{n+2}
\dd\tau
\notag
\\
&\quad
+
s2^{s-1}
\sum_{n=0}^{\infty}
|c_n|(2\pi)^{-(n+1)}
\int_{T_1}^{T_2}
\int_{\xi_1+\cdots+\xi_{n+3}=0}
|\xi_1^{(1)}|
|\xi_1|^{s-1}
|\widehat\theta(\xi_1)|
\notag
\\*[-2pt]
&\qquad
\times
|\xi_2^{(1)}|
|\widehat\theta(\xi_2)|
|\xi_3|^s
|\widehat\theta(\xi_3)|
\prod_{k=4}^{n+3}
|\xi_k|
|\widehat\theta(\xi_k)|
\dd\xi_1
\cdots
\dd\xi_{n+2}
\dd\tau
\notag
\\
&
\le
2s2^{s-1}
\sum_{n=0}^{\infty}
|c_n|(2\pi)^{-(n+1)}
\int_{T_1}^{T_2}
\int_{\xi_1+\cdots+\xi_{n+3}=0}
|\xi_1^{(1)}|
|\widehat\theta(\xi_1)|
\notag
\\*[-2pt]
&\qquad
\times
|\xi_2^{(1)}|
|\xi_2|^{s-1}
|\widehat\theta(\xi_2)|
|\xi_3|^s
|\widehat\theta(\xi_3)|
\prod_{k=4}^{n+3}
|\xi_k|
|\widehat\theta(\xi_k)|
\dd\xi_1
\cdots
\dd\xi_{n+2}
\dd\tau.
\label{eq:I6start}
\end{align}

From \eqref{eq:W1bound},
\begin{equation}
\label{eq:I6integral}
\begin{aligned}
&
\int_{\xi_1+\cdots+\xi_{n+3}=0}
|\xi_1^{(1)}|
|\widehat\theta(\xi_1)|
|\xi_2^{(1)}|
|\xi_2|^{s-1}
|\widehat\theta(\xi_2)|
|\xi_3|^s
|\widehat\theta(\xi_3)|
\\[-2pt]
&\qquad
\times
\prod_{k=4}^{n+3}
|\xi_k|
|\widehat\theta(\xi_k)|
\dd\xi_1
\cdots
\dd\xi_{n+2}
\\
&\quad
\le
2(n+2)
\left(
\sqrt{\frac{\pi}{s-2}}
\norm{\theta}_{H^s}
\right)^{n+1}
\norm{\mathcal R_1\theta}_{H^s}^2.
\end{aligned}
\end{equation}
Hence, we obtain
\begin{align}
I_6(T_1,T_2)
&
\le
2s2^{s-1}
\sum_{n=0}^{\infty}
|c_n|(n+2)^2
\notag
\\*
&\qquad
\times
\left(
\frac{\norm{\theta}_{L^\infty(T_1,T_2;H^s)}}{2\sqrt{\pi(s-2)}}
\right)^{n+1}
\norm{\mathcal R_1\theta}_{L^2(T_1,T_2;H^s)}^2
\notag
\\
&
=
s2^s
\sum_{n=1}^{\infty}
|c_{n-1}|(n+1)^2
\notag
\\*
&\qquad
\times
\left(
\frac{\norm{\theta}_{L^\infty(T_1,T_2;H^s)}}{2\sqrt{\pi(s-2)}}
\right)^n
\norm{\mathcal R_1\theta}_{L^2(T_1,T_2;H^s)}^2
\notag
\\
&
\le
s2^s
\sum_{n=1}^{\infty}
|c_{n-1}|(n+1)^2
\alpha^n
\norm{\mathcal R_1\theta}_{L^2(T_1,T_2;H^s)}^2.
\label{eq:I6}
\end{align}
\paragraph{Collection of the estimates}
Adding \eqref{eq:I1}, \eqref{eq:I2}, \eqref{eq:I3}, \eqref{eq:I4}, \eqref{eq:I5}, and \eqref{eq:I6}, and then applying \eqref{eq:cn-bound}, we obtain
\begin{equation}
\label{eq:recorded-close}
\begin{aligned}
\sum_{j=1}^6
I_j(T_1,T_2)
\le{}
&
2s(1-\alpha)^{-(2s+1)}
\alpha^2
\norm{\mathcal R_1\theta}_{L^2(T_1,T_2;H^s)}^2
\\
&
+
s2^{s+1}
\alpha
\norm{\mathcal R_1\theta}_{L^2(T_1,T_2;H^s)}^2
\\
&
+
3
\sum_{n=1}^{\infty}
|c_n|(n+1)^2
\alpha^n
\norm{\mathcal R_1\theta}_{L^2(T_1,T_2;H^s)}^2
\\
&
+
11s(s-1)
2^{s-2}
\sum_{n=1}^{\infty}
|c_{n-1}|(n+1)^2
\alpha^n
\norm{\mathcal R_1\theta}_{L^2(T_1,T_2;H^s)}^2
\\
&
+
s2^s
\sum_{n=1}^{\infty}
|c_{n-1}|(n+1)^2
\alpha^n
\norm{\mathcal R_1\theta}_{L^2(T_1,T_2;H^s)}^2
\\
\le{}
&
2s(1-\alpha)^{-(2s+1)}
\alpha^2
\norm{\mathcal R_1\theta}_{L^2(T_1,T_2;H^s)}^2
\\
&
+
4s^2
2^s
\sum_{n=1}^{\infty}(|c_n|+|c_{n-1}|)
(n+1)^2
\alpha^n
\norm{\mathcal R_1\theta}_{L^2(T_1,T_2;H^s)}^2
\\
\leq{}
&
\Big[
2s(1-\alpha)^{-(2s+1)}
\alpha^2
+
s^2
2^{s+2}
C^*_{\alpha, s}
\Big]
\\
&
\times
\|\mathcal R_1\theta\|_{L^2(T_1,T_2;H^s)}^2.
\end{aligned}
\end{equation}
Here, we used $c_0=1$, $c_1=2s$, and $s>2$ to absorb the $I_2$ and $I_3$ bounds into the $n=1$ summand. Inserting \eqref{eq:recorded-close} into
\eqref{eq:Ds-source-bound} and
\eqref{eq:Ps-reverse-source-bound} respectively, we obtain
\[
\begin{aligned}
D_s(T_1,T_2)^2
\leq{}
&
\frac32
P_s(T_1)^2
+
\Big[
6s(1-\alpha)^{-(2s+1)}
\alpha^2
+
3s^2
2^{s+2}
C^*_{\alpha, s}
\Big]
\times
\|\mathcal R_1\theta\|_{L^2(T_1,T_2;H^s)}^2,
\end{aligned}
\]
\[
\begin{aligned}
P_s(T_1)^2
\leq{}
&
P_s(T_2)^2
+
\Big[
2
+
4s(1-\alpha)^{-(2s+1)}
\alpha^2
+
s^2
2^{s+3}
C^*_{\alpha, s}
\Big]
\times
\|\mathcal R_1\theta\|_{L^2(T_1,T_2;H^s)}^2.
\end{aligned}
\]
\end{proof}

\section{\texorpdfstring{$L^2$}{L2}-convergence lemma and density coordinates}

In this section, we establish an $L^2$-convergence lemma
and apply it to equation~\eqref{equ1.2}
via a density-coordinate transformation.

\subsection{\texorpdfstring{$L^2$}{L2}-convergence lemma}

The time-integrated dissipation contained in $D_s$
motivates the following convergence lemma.

\begin{Lemma}\label{lem:l2-convergence}
Let $z\in C([0,\infty);L^2(\mathbb R^2))$ be a
solution of
\begin{equation}
\partial_t z
-
\mathcal R_1^2
z
=
\mathcal R_1
f,
\label{eq:l2-convergence-equation}
\end{equation}
where
\[
f
\in
L^2(0,\infty;L^2(\mathbb R^2)).
\]
Then
\[
\lim_{t\to\infty}
\|z(t)\|_{L^2}
=
0.
\]
\end{Lemma}

\begin{proof}
If $h,h'\in L^1(0,\infty)$, then
\begin{equation}
\lim_{t\to\infty}
h(t)
=
0.
\label{eq:l1-convergence-fact}
\end{equation}

Multiplying \eqref{eq:l2-convergence-equation} by
$z$ and integrating over $\mathbb R^2$, we obtain
\[
\begin{aligned}
\frac12
\frac{d}{dt}
\|z\|_{L^2}^2
+
\|\mathcal R_1z\|_{L^2}^2
&
=
\int_{\mathbb R^2}
\mathcal R_1
f\,
z\,
dx
\\
&
=
-
\int_{\mathbb R^2}
f\,
\mathcal R_1
z\,
dx
\\
&
\leq
\frac12
\|f\|_{L^2}^2
+
\frac12
\|\mathcal R_1z\|_{L^2}^2.
\end{aligned}
\]
Hence,
\[
\frac{d}{dt}
\|z\|_{L^2}^2
+
\|\mathcal R_1z\|_{L^2}^2
\leq
\|f\|_{L^2}^2.
\]
Integrating over $[0,t]$, we obtain
\[
\begin{aligned}
\|z(t)\|_{L^2}^2
&
+
\int_0^t
\|\mathcal R_1z(\tau)\|_{L^2}^2\,
d
\tau
\leq
\|z(0)\|_{L^2}^2
+
\int_0^t
\|f(\tau)\|_{L^2}^2\,
d
\tau.
\end{aligned}
\]

For $0<\delta<1$, we define the Fourier projections
\[
P_\delta
f
=
\mathcal F^{-1}
\left(
\mathbf 1_\delta
\left(
\frac{\xi^{(1)}}{|\xi|}
\right)
\widehat f(\xi)
\right),
\qquad
P_\delta^\perp
f
=
f
-
P_\delta
f,
\]
where
\[
\mathbf 1_\delta(a)
=
\begin{cases}
1,
&
|a|
\leq
\delta,
\\
0,
&
|a|
>
\delta.
\end{cases}
\]
Applying $P_\delta^\perp$ to
\eqref{eq:l2-convergence-equation}, we obtain
\[
\frac{d}{dt}
\|P_\delta^\perp z\|_{L^2}^2
+
\|\mathcal R_1P_\delta^\perp z\|_{L^2}^2
\leq
\|f\|_{L^2}^2.
\]
Since
\[
\|\mathcal R_1P_\delta^\perp z\|_{L^2}^2
\geq
\delta^2
\|P_\delta^\perp z\|_{L^2}^2,
\]
it follows that
\[
\frac{d}{dt}
\|P_\delta^\perp z\|_{L^2}^2
+
\delta^2
\|P_\delta^\perp z\|_{L^2}^2
\leq
\|f\|_{L^2}^2.
\]
Hence,
\[
\|P_\delta^\perp z(t)\|_{L^2}^2,
\quad
\frac{d}{dt}
\|P_\delta^\perp z(t)\|_{L^2}^2
\in
L^1(0,\infty).
\]
Applying \eqref{eq:l1-convergence-fact},
we conclude that
\begin{equation}
\lim_{t\to\infty}
\|P_\delta^\perp z(t)\|_{L^2}^2
=
0.
\label{eq:l2-high-frequency-convergence}
\end{equation}

Next, applying $P_\delta$ to
\eqref{eq:l2-convergence-equation}
and integrating the corresponding energy
inequality over $[T,t]$, we obtain
\[
\begin{aligned}
\|P_\delta z(t)\|_{L^2}^2
&
+
\int_T^t
\|\mathcal R_1P_\delta z(\tau)\|_{L^2}^2\,
d
\tau
\leq
\|P_\delta z(T)\|_{L^2}^2
+
\int_T^t
\|f(\tau)\|_{L^2}^2\,
d
\tau.
\end{aligned}
\]
Consequently, for every $T\geq0$,
\begin{equation}
\begin{aligned}
\limsup_{t\to\infty}
\|P_\delta z(t)\|_{L^2}^2
\leq{}
&
\|P_\delta z(T)\|_{L^2}^2
+
\int_T^\infty
\|f(\tau)\|_{L^2}^2\,
d
\tau.
\end{aligned}
\label{eq:l2-low-frequency-estimate}
\end{equation}
It is clear that
\[
\|z(t)\|_{L^2}^2
=
\|P_\delta z(t)\|_{L^2}^2
+
\|P_\delta^\perp z(t)\|_{L^2}^2.
\]
Combining this identity with
\eqref{eq:l2-high-frequency-convergence}
and \eqref{eq:l2-low-frequency-estimate},
we obtain
\[
\limsup_{t\to\infty}
\|z(t)\|_{L^2}^2
\leq{}
\|P_\delta z(T)\|_{L^2}^2
+
\int_T^\infty
\|f(\tau)\|_{L^2}^2\,
d
\tau.
\]
Letting $\delta\to0$, we have
\[
\|P_\delta z(T)\|_{L^2}^2
\longrightarrow
0
\]
by the dominated convergence theorem, since the set
\[
\{
\xi
\in
\mathbb R^2
:
\xi^{(1)}
=
0
\}
\]
has zero two-dimensional Lebesgue measure.
Thus,
\[
\limsup_{t\to\infty}
\|z(t)\|_{L^2}^2
\leq
\int_T^\infty
\|f(\tau)\|_{L^2}^2\,
d
\tau.
\]
Finally, letting $T\to\infty$, we obtain
\[
\lim_{t\to\infty}
\|z(t)\|_{L^2}
=
0.
\]
\end{proof}

\subsection{Density coordinates}

To apply Lemma~3.1, we introduce the
density-coordinate transformation $y=\Phi(t,x)$,
defined by
\[
y_1
=
x_1,
\qquad
y_2
=
-
x_2
+
\theta(t,x).
\]
If $\|\partial_2\theta\|_{L^\infty}<1$
and $\|\theta\|_{L^\infty}<\infty$, then
$\Phi(t,\cdot)$ is invertible, and its inverse
satisfies
\[
\Phi^{-1}(t,y)
=
\bigl(
y_1,
\theta(t,\Phi^{-1}(t,y))
-
y_2
\bigr).
\]
Direct differentiation gives
\[
D_y
\Phi^{-1}
=
\begin{pmatrix}
1
&
0
\\[2mm]
-
\dfrac{\partial_1\theta}
{-1+\partial_2\theta}
&
\dfrac{1}{-1+\partial_2\theta}
\end{pmatrix}
\circ
\Phi^{-1},
\]
and
\[
\partial_t \Phi^{-1}
=
\left(
0,
-
\frac{\partial_t\theta}
{-1+\partial_2\theta}
\right)
\circ
\Phi^{-1}.
\]
Let $\psi$ be a stream function satisfying
\begin{equation}
-
\Delta
\psi
=
\partial_1 \theta.
\label{eq:density-stream}
\end{equation}
We define
\[
\widetilde\theta
=
\theta
\circ
\Phi^{-1},
\qquad
\widetilde\psi
=
\psi
\circ
\Phi^{-1}.
\]

Recalling that
\[
u
=
(-\mathcal R_1\mathcal R_2\theta,\mathcal R_1^2\theta)
=
(-\partial_2\psi,\partial_1\psi),
\]
the chain rule and equation~\eqref{equ1.2} yield
\begin{equation}
\begin{aligned}
\partial_t \widetilde\theta
&
=
(\partial_t\theta)
\circ
\Phi^{-1}
+
(\nabla\theta)
\circ
\Phi^{-1}
\cdot
\partial_t \Phi^{-1}
\\
&
=
\left(
\frac{-\partial_t\theta}
{-1+\partial_2\theta}
\right)
\circ
\Phi^{-1}
\\
&
=
\left(
u^2
+
\frac{u^1\partial_1\theta}
{-1+\partial_2\theta}
\right)
\circ
\Phi^{-1}
\\
&
=
\partial_{y_1} \widetilde\psi.
\end{aligned}
\label{eq:density-time}
\end{equation}
Next, let $\chi\in C_c^\infty(\mathbb R^2)$.
Multiplying \eqref{eq:density-stream} by
$\chi\circ\Phi$ and integrating by parts,
we obtain
\begin{equation}
\begin{aligned}
\int_{\mathbb R^2}
\partial_1 \theta\,
(\chi\circ\Phi)\,
dx
&
=
-
\int_{\mathbb R^2}
\Delta
\psi\,
(\chi\circ\Phi)\,
dx
\\
&
=
\int_{\mathbb R^2}
\nabla \psi
\cdot
\nabla(\chi\circ\Phi)\,
dx.
\end{aligned}
\label{eq:density-weak-original}
\end{equation}
If $\|\partial_2\theta\|_{L^\infty}<1$,
we have
\[
(\partial_1\theta)
\circ
\Phi^{-1}\,
|\det D_y\Phi^{-1}|
=
\partial_{y_1} \widetilde\theta.
\]
Making the change of variables
$x=\Phi^{-1}(t,y)$ in
\eqref{eq:density-weak-original},
we obtain
\[
\begin{aligned}
\int_{\mathbb R^2}
\partial_{y_1} \widetilde\theta\,
\chi\,
dy
&
=
\int_{\mathbb R^2}
\nabla_y \widetilde\psi
\cdot
\nabla_y \chi\,
dy
+
\int_{\mathbb R^2}
(\nabla_y\widetilde\psi\,B)
\cdot
\nabla_y \chi\,
dy,
\end{aligned}
\]
where
\[
\begin{aligned}
B
&
=
\left[
\frac{D_x\Phi(D_x\Phi)^T}
{|\det D_x\Phi|}
\right]
\circ
\Phi^{-1}
-
I
\\[2mm]
&
=
\left[
\frac{1}{1-\partial_2\theta}
\begin{pmatrix}
\partial_2 \theta
&
\partial_1 \theta
\\[2mm]
\partial_1 \theta
&
(\partial_1\theta)^2
-
\partial_2 \theta
+
(\partial_2\theta)^2
\end{pmatrix}
\right]
\circ
\Phi^{-1}.
\end{aligned}
\]
Thus, we obtain
\begin{equation}
-
\Delta_y
\widetilde\psi
=
\partial_{y_1} \widetilde\theta
+
\operatorname{div}_y
(\nabla_y\widetilde\psi\,B).
\label{eq:density-elliptic}
\end{equation}
Substituting \eqref{eq:density-elliptic}
into \eqref{eq:density-time}, we obtain
\begin{equation}
\partial_t \widetilde\theta
=
\mathcal R_{1,y}^2
\widetilde\theta
+
\mathcal R_{1,y}(-\Delta_y)^{-1/2}
\operatorname{div}_y
(\nabla_y\widetilde\psi\,B).
\label{eq:density-evolution}
\end{equation}
where
\[
\mathcal R_{1,y}
:=
\partial_{y_1}(-\Delta_y)^{-1/2}.
\]

We now establish the following convergence result.
\begin{Proposition}
\label{prop:density-convergence}
Let $s>2$ and $0<\delta<1$.
Suppose that
\[
\theta
\in
C([0,\infty);
H^s(\mathbb R^2)
)
\]
is a global solution of \eqref{equ1.2}
satisfying
\begin{equation}
\begin{gathered}
\|\partial_2\theta\|_
{L^\infty(0,\infty;L^\infty)}
\leq
\delta,
\\
\|\theta\|_
{L^\infty(0,\infty;H^s)}
<
\infty,
\\
\|\mathcal R_1\theta\|_
{L^2(0,\infty;L^2)}
<
\infty.
\end{gathered}
\label{eq:density-assumptions}
\end{equation}
Then, for every $0\leq s'<s$,
\begin{equation}
\lim_{t\to\infty}
\|\theta(t)\|_{H^{s'}}
=
0.
\label{eq:density-convergence}
\end{equation}
\end{Proposition}

\begin{proof}
By the change of variables $y=\Phi(t,x)$,
we have
\begin{equation}
\begin{aligned}
\|\widetilde\theta(t)\|_{L_y^2}^2
&
=
\int_{\mathbb R^2}
(\theta\circ\Phi^{-1})^2\,
dy
\\
&
=
\int_{\mathbb R^2}
\theta^2(1-\partial_2\theta)\,
dx
\\
&
=
\int_{\mathbb R^2}
\theta^2\,
dx
-
\frac13
\int_{\mathbb R^2}
\partial_2(\theta^3)\,
dx
\\
&
=
\int_{\mathbb R^2}
\theta^2\,
dx,
\end{aligned}
\label{eq:density-L2-identity}
\end{equation}
where, in the last equality, we used
$\theta^3\in W^{1,1}(\mathbb R^2)$,
which follows from
$\theta\in H^s(\mathbb R^2)$ with $s>2$.

By \eqref{eq:density-evolution},
the function $\widetilde\theta$
satisfies an equation of the form
\eqref{eq:l2-convergence-equation}, with
\[
f
=
(-\Delta_y)^{-1/2}
\operatorname{div}_y
(\nabla_y\widetilde\psi\,B).
\]
Under the assumption \eqref{eq:density-assumptions},
the transformed coefficients are uniformly
bounded, and
\[
\begin{aligned}
\|f(t)\|_{L_y^2}
\leq
\|\nabla_y\widetilde\psi\,B\|_{L_y^2}
\leq
C
\|\nabla_y\widetilde\psi\|_{L_y^2}
\leq
C
\|u(t)\|_{L_x^2}
=
C
\|\mathcal R_1\theta(t)\|_{L_x^2}.
\end{aligned}
\]
Consequently,
\[
f
\in
L^2(0,\infty;L^2(\mathbb R^2)).
\]
Applying Lemma~3.1 to
\eqref{eq:density-evolution}, we obtain
\[
\lim_{t\to\infty}
\|\widetilde\theta(t)\|_{L_y^2}
=
0.
\]
It follows from
\eqref{eq:density-L2-identity} that
\[
\lim_{t\to\infty}
\|\theta(t)\|_{L_x^2}
=
0.
\]
Finally, for every $0\leq s'<s$,
the interpolation inequality and
\eqref{eq:density-assumptions} give
\[
\begin{aligned}
\|\theta(t)\|_{H^{s'}}
&
\leq
\|\theta(t)\|_{L^2}^{1-s'/s}
\|\theta(t)\|_{H^s}^{s'/s}
\longrightarrow
0
\qquad
\text{as }
t
\to
\infty.
\end{aligned}
\]
\end{proof}

\section{Proof of Theorem~\ref{thm:main}}
\label{sec:proof-main}

\begin{proof}
Let $0<\alpha_0<1$ be sufficiently small such that
\[
C_{\alpha_0,s}^2
\Big[
6s(1-\alpha_0)^{-(2s+1)}
\alpha_0^2
+
3s^2
2^{s+2}
C^*_{\alpha_0, s}
\Big]
\leq
\frac14,
\]
and define
\[
\epsilon_s
=
\frac{\alpha_0}{C_{\alpha_0,s}^{2}}
\sqrt{\pi\left(\frac{s+2}{2}-2\right)},
\]
where the constants $C_{\alpha,s}$ and $C^*_{\alpha, s}$ are defined in
\eqref{eq:energy-equivalence-constant} and \eqref{new1}, respectively.

Let $T^*\in(0,\infty]$ be the maximal existence time
of the solution to \eqref{equ1.2}, with
\[
\theta
\in
C([0,T];H^s),
\qquad
\mathcal R_1
\theta
\in
L^2(0,T;H^s)
\]
for every $T<T^*$. We first prove that $T^*=\infty$ whenever
$\|\theta_0\|_{H^s}\leq\ep_s$. To this end, we define
\[
T^{**}
=
\sup
\left\{
0
<
T
<
T^*
:
E_s(T)
\leq
2
\alpha_0
\sqrt{\pi\left(\frac{s+2}{2}-2\right)}
\right\}.
\]
For $0<T<T^{**}$, we have
\[
\begin{aligned}
\|\theta\|_{L^\infty(0,T;H^s)}
\leq
E_s(T)
\leq
2
\alpha_0
\sqrt{\pi\left(\frac{s+2}{2}-2\right)}.
\end{aligned}
\]
By Proposition~\ref{Pro1}
and \eqref{eq:Ps-equivalence}--\eqref{eq:Ds-Es-equivalence},
\[
\begin{aligned}
E_s(T)^2
\leq{}
&
\frac32
C_{\alpha_0,s}^4
\|\theta_0\|_{H^s}^2
+
C_{\alpha_0,s}^2
\Big[
6s(1-\alpha_0)^{-(2s+1)}
\alpha_0^2
+
3s^2
2^{s+2}
C^*_{\alpha_0, s}
\Big]
E_s(T)^2
\\
\leq{}
&
\frac32
C_{\alpha_0,s}^4
\|\theta_0\|_{H^s}^2
+
\frac14
E_s(T)^2.
\end{aligned}
\]
Consequently,
\[
\begin{aligned}
E_s(T)^2
&
\leq
2C_{\alpha_0,s}^4
\|\theta_0\|_{H^s}^2
\leq
2C_{\alpha_0,s}^4
\epsilon_s^2
\\
&
=
2
\alpha_0^2
\pi
\left(
\frac{s+2}{2}
-
2
\right)
\\
&
<
4
\alpha_0^2
\pi
\left(
\frac{s+2}{2}
-
2
\right).
\end{aligned}
\]
It follows by continuity that $T^{**}=T^*$.
Moreover, if $T^*<\infty$, the uniform $H^s$ bound implies
\[
\int_0^{T^*}
\left(
\|\nabla\theta(t)\|_{L^\infty}
+
\|\nabla u(t)\|_{L^\infty}
\right)\,
dt
<
\infty.
\]
The local solution can therefore be extended beyond
$T^*$, yielding a contradiction. Thus we obtain $T^*=\infty$.

We have proved the global estimate
\begin{equation}
E_s(0,\infty)
\leq
2
\alpha_0
\sqrt{\pi\left(\frac{s+2}{2}-2\right)}.
\label{eq:global-energy-smallness}
\end{equation}
Next, we establish strong convergence.
For every
\[
\frac{s+2}{2}
\leq
s
'
<
s,
\]
it follows from \eqref{eq:global-energy-smallness} that
\begin{equation}
\begin{aligned}
\|\theta\|_{L^\infty(0,\infty;H^{s'})}
&
\leq
E_{s'}(0,\infty)
\leq
E_s(0,\infty)
\\
&
\leq
2
\alpha_0
\sqrt{\pi\left(\frac{s+2}{2}-2\right)}
\\
&
\leq
2
\alpha_0
\sqrt{\pi(s'-2)}.
\end{aligned}
\label{eq:lower-order-smallness}
\end{equation}
By \eqref{eq:L1theta},
\begin{equation}\label{eq:global-gradient-smallness}
\|\partial_2\theta\|_{L^\infty(0,\infty;L^\infty)}
\leq
\frac{\|\theta\|_{L^\infty(0,\infty;H^{s'})}}
{2\sqrt{\pi(s'-2)}}
\leq
\alpha_0.
\end{equation}
From \eqref{eq:Ps-equivalence},
Proposition~\ref{prop:density-convergence},
\eqref{eq:global-energy-smallness}
and \eqref{eq:global-gradient-smallness},
we obtain
\begin{equation}
\lim_{t\to\infty}
P_r(t)^2
=
0,
\qquad
0
\leq
r
<
s.
\label{eq:lower-weighted-convergence}
\end{equation}

Now, let $\frac{s+2}{2}\leq s'<s$ and $0\leq T<t$.
By Proposition~\ref{Pro1},
\eqref{eq:global-energy-smallness},
and \eqref{eq:global-gradient-smallness},
\[
\begin{aligned}
P_{s'}(T)^2
\leq{}
&
P_{s'}(t)^2
+
\Big[
2
+
4s
'(1-\alpha_0)^{-(2s'+1)}
\alpha_0^2
+
(s')^2
2^{s'+3}
C^*_{\alpha_0, s'}
\Big]
\|\mathcal R_1\theta\|_{L^2(T,t;H^{s'})}^2
\\
\leq{}
&
P_{s'}(t)^2
+
\Big[
2
+
4s(1-\alpha_0)^{-(2s+1)}
\alpha_0^2
+
s^2
2^{s+3}
C^*_{\alpha_0, s}
\Big]
\|\mathcal R_1\theta\|_{L^2(T,t;H^s)}^2.
\end{aligned}
\]
Sending $t\to\infty$ and using
\eqref{eq:lower-weighted-convergence}, we obtain
\begin{align}\label{new2}
P_{s'}(T)^2
\leq
\Big[
2
+
4s(1-\alpha_0)^{-(2s+1)}
\alpha_0^2
+
s^2
2^{s+3}
C^*_{\alpha_0, s}
\Big]
\|\mathcal R_1\theta\|_{L^2(T,\infty;H^s)}^2.
\end{align}
For each fixed $T\geq0$, since $\theta(T)\in H^s$,
the dominated convergence theorem yields
\[
\lim_{s'\uparrow s}
\|\Lambda_i^{s'}\theta(T)-\Lambda_i^s\theta(T)\|_{L^2}
=0,
\qquad i=1,2.
\]
By the mean value theorem and
\eqref{eq:global-gradient-smallness},
\[
\|\theta_{s'}(T)-\theta_s(T)\|_{L^\infty}
\leq C|s'-s|
\longrightarrow0
\qquad\text{as }s'\uparrow s,
\]
where $C$ depends only on $s$ and $\alpha_0$.
Hence,
\[
\begin{aligned}
&\left|
\int_{\mathbb R^2}
(\Lambda_2^{s'}\theta(T))^2\theta_{s'}(T)\,dx
-
\int_{\mathbb R^2}
(\Lambda_2^s\theta(T))^2\theta_s(T)\,dx
\right|\\
&\quad\leq
\|\theta_{s'}(T)-\theta_s(T)\|_{L^\infty}
\|\Lambda_2^{s'}\theta(T)\|_{L^2}^2\\
&\qquad+
\|\theta_s(T)\|_{L^\infty}
\|\Lambda_2^{s'}\theta(T)-\Lambda_2^s\theta(T)\|_{L^2}\\
&\qquad\quad\times
\Big(
\|\Lambda_2^{s'}\theta(T)\|_{L^2}
+\|\Lambda_2^s\theta(T)\|_{L^2}
\Big)\\
&\quad\longrightarrow0
\qquad\text{as }s'\uparrow s.
\end{aligned}
\]
These convergences imply
\[
\lim_{s'\uparrow s}P_{s'}(T)^2=P_s(T)^2.
\]
Letting $s' \to  s$ in \eqref{new2} gives
\[
P_s(T)^2
\leq
\Big[
2
+
4s(1-\alpha_0)^{-(2s+1)}
\alpha_0^2
+
s^2
2^{s+3}
C^*_{\alpha_0, s}
\Big]
\|\mathcal R_1\theta\|_{L^2(T,\infty;H^s)}^2.
\]
Since
\[
\mathcal R_1
\theta
\in
L^2(0,\infty;H^s),
\]
letting $T\to\infty$ gives
\[
\lim_{T\to\infty}
P_s(T)^2
=
0.
\]
Finally, using \eqref{eq:Ps-equivalence}, we conclude that
\[
\lim_{T\to\infty}
\|\theta(T)\|_{H^s}
=
0.
\]
%The corresponding convergence of $u$ follows from its
%representation in \eqref{equ1.2}. 
This completes the proof.
\end{proof}

\bigskip

\noindent{\bf{Acknowledgement.}} The work of Y.L. was supported by National Natural Science Foundation of China (12571228), Natural Science Foundation of Anhui Province (2408085MA018). The work of S.L. was supported by National Natural Science Foundation of China (12601425),  Tianyuan Fund for Mathematics of the National Natural Science Foundation of China (12526541).

\vspace{4mm}
\noindent{\bf{Data Availability.}} Data sharing is not applicable to this article as no datasets were generated or analyzed
during the current study. 
\vspace{4mm}

\section*{Declarations}

\noindent{\bf{Conflicts of interest.}} All authors certify that there are no conflicts of interest for this work.

\vspace{4mm}

\noindent{\bf{Ethics approval.}} This article does not contain any studies involving humans or animals.

\end{document}